# Asymptotic expansion and central limit theorem for quadratic variations of Gaussian processes

ARNAUD BEGYN

*Lycée Pierre de Fermat, 31000 Toulouse, France. E-mail: arnaud.begyn@free.fr*

Cohen, Guyon, Perrin and Pontier have given assumptions under which the second-order quadratic variations of a Gaussian process converge almost surely to a deterministic limit. In this paper we present two new convergence results about these variations: the first is a deterministic asymptotic expansion; the second is a central limit theorem. Next we apply these results to identify two-parameter fractional Brownian motion and anisotropic fractional Brownian motion.

## Introduction

In this paper we consider second-order quadratic variations of a Gaussian process $X$. From Cohen *et al.* [10] we know that it converges to a deterministic limit under convenient conditions on the covariance function of the process. First we sharpen this result: we show that if an asymptotic expansion of the covariance function is known, we get an asymptotic expansion of the second-order quadratic variation. Next we establish a central limit theorem related to the previous result. We apply these results to two-parameter fractional Brownian motion, which is a generalization of fractional Brownian motion that has non-stationary increments, and to anisotropic fractional Brownian field, which is a multidimensional anisotropic generalization of fractional Brownian motion.

In the first section, we state the notation. In the second section, we prove the main theorems about the second-order quadratic variation. In the third section, we study the case of two-parameter fractional Brownian motion. In the fourth section, we consider anisotropic fractional Brownian motion.

# 1. Notation

Let $X = \{X_t; t \in [0,1]\}$ be a Gaussian process. We denote by $t \mapsto M_t$ its mean function and by $(s,t) \mapsto R(s,t)$ its covariance function.

We define the second-order increments of $R$ as

$$\delta_1^h R(s,t) = R(s+h,t) + R(s-h,t) - 2R(s,t),$$

$$\delta_2^h R(s,t) = R(s,t+h) + R(s,t-h) - 2R(s,t).$$

For fractional processes (i.e., processes whose properties are close to those of fractional Brownian motion), we use the second-order quadratic variation

$$V_n(X) = \sum_{k=1}^{n-1} [X_{(k+1)/n} + X_{(k-1)/n} - 2X_{k/n}]^2, \tag{1}$$

because the standard quadratic variation does not satisfy a central limit theorem in general.

To be sure that $V_n(X)$ converges almost surely to a deterministic limit, we need to normalize this quantity. A result of the form

$$\lim_{n \to +\infty} n^{1-\gamma} V_n(X) = \int_0^1 g(t)\, dt \qquad \text{a.s.} \tag{2}$$

is expected, where $\gamma$ is related to the regularity of the paths of $X$, and $g$ is related to the non-differentiability of $R$ on the diagonal $\{s = t\}$ and is called the singularity function of the process. In this paper, we consider a class of processes for which a more general normalization is needed. Moreover, we present a better result because we give an asymptotic expansion of the left-hand side of (2).

We will say that a Borel function $\psi : ]0,a[ \to \mathbb{R}$ $(a > 0)$ is regularly varying with index $\beta \in \mathbb{R}$ if $\psi(h) = h^\beta L(h)$, where $L$ is a slowly varying function

$$\forall \lambda > 0 \qquad \lim_{x \to 0^+} \frac{L(\lambda x)}{L(x)} = 1.$$

Let $d \in \mathbb{N}^*$. Standard fractional Brownian motion (FBM) $B^H = \{B_t^H; t \in \mathbb{R}^d\}$, with Hurst index $H \in ]0,1[$, is the unique continuous centered Gaussian process, which has the covariance function

$$\forall s, t \in \mathbb{R}^d \qquad \text{Cov}(B_s^H, B_t^H) = \tfrac{1}{2}(|s|^{2H} + |t|^{2H} - |s-t|^{2H}), \tag{3}$$

where $|\cdot|$ denotes the Euclidean norm.

In next section, we use the notation (we drop the superscript index $n$ wherever it is possible)

$$\Delta X_k^{(n)} = X_{(k+1)/n} + X_{(k-1)/n} - 2X_{k/n}, \qquad k = 1, \ldots, n-1,$$

and

$$d_{jk}^{(n)} = \mathbb{E}(\Delta X_j^{(n)} \Delta X_k^{(n)}), \qquad j, k = 1, \ldots, n-1. \tag{4}$$

## 2. The results

In this section we sharpen (2). First, we prove a deterministic asymptotic expansion of $V_n(X)$ under certain conditions on the covariance function. Second, we prove a central limit theorem.

Examples of the application of Theorem 2.1 with a non-trivial slowly varying function $L(h)$ can be found in Section 4.2.

### 2.1. Asymptotic expansion

**Theorem 1.** *Assume that $X$ satisfies the following statements:*

1. *$t \mapsto M_t = \mathbb{E}X_t$ has a bounded first derivative in $[0, 1]$.*

2. *The covariance function $R$ has the following properties:*

   (a) *$R$ is continuous in $[0, 1]^2$.*

   (b) *The derivative $\frac{\partial^4 R}{\partial s^2 \partial t^2}$ exists and is continuous in $]0, 1]^2 \setminus \{s = t\}$. There exist a constant $C > 0$, a real $\gamma \in {]}0, 2{[}$ and a positive slowly varying function $L : {]}0, 1{[} \to {]}0, +\infty{[}$ such that*

   $$\forall s, t \in {]}0, 1]^2 \setminus \{s = t\} \qquad \left| \frac{\partial^4 R}{\partial s^2 \partial t^2}(s, t) \right| \leq C \frac{L(|s-t|)}{|s-t|^{2+\gamma}}. \tag{5}$$

   (c) *There exist $q+1$ functions ($q \in \mathbb{N}$) $g_0, g_1, \ldots, g_q$ from $]0, 1[$ to $\mathbb{R}$, $q$ real numbers $0 < \nu_1 < \cdots < \nu_q$ and a function $\phi : {]}0, 1{[} \to {]}0, +\infty{[}$ such that:*

      (i) *if $q \geq 1$, then $\forall 0 \leq i \leq q-1$, $g_i$ is Lipschitz on $]0, 1[$;*

      (ii) *$g_q$ is bounded on $]0, 1[$;*

      (iii) *we have*

      $$\sup_{h \leq t \leq 1-h} \left| \frac{(\delta_1^h \circ \delta_2^h R)(t, t)}{h^{2-\gamma} L(h)} - g_0(t) - \sum_{i=1}^q g_i(t) \phi(h)^{\nu_i} \right| \overset{h \to 0^+}{=} o(\phi(h)^{\nu_q}), \tag{6}$$

      *where if $q = 0$, then $\sum_{i=1}^q g_i(t)\phi(h)^{\nu_i} = 0$ and $\phi(h)^{\nu_q} = 1$; else if $q \neq 0$, then $\lim_{h \to 0^+} \phi(h) = 0$.*

3. *If $q \neq 0$, we assume that*

   $$\lim_{n \to +\infty} \frac{\log n}{n \phi(1/n)^{\nu_q}} = 0. \tag{7}$$

4. *If $X$ is not centered, we make the additional assumption*

   $$\lim_{n \to +\infty} \frac{1}{n^\gamma L(1/n) \phi(1/n)^{\nu_q}} = 0, \tag{8}$$

*where if $q = 0$, then $\phi(1/n)^{\nu_q} = 1$.*

*Then, for all $t \in [0, 1]$, we have almost surely*

$$\frac{n^{1-\gamma}}{L(1/n)} V_n(X) \overset{n \to +\infty}{=} \int_0^1 g_0(x)\,dx + \sum_{i=1}^q \left(\int_0^1 g_i(x)\,dx\right)\phi\left(\frac{1}{n}\right)^{\nu_i} + o\left(\phi\left(\frac{1}{n}\right)^{\nu_q}\right). \quad (9)$$

**Remarks.** (i) If the assumption (6) is fulfilled for $q_*$, then it is fulfilled for all $q \in \{0, 1, \ldots, q_*\}$ too with the truncated sequences $(g_i)_{0 \le i \le q}$ and $(\nu_i)_{0 \le i \le q}$. The maximal value of $q$ is given by the assumption (8), which yields an upper bound for the value of $\nu_q$.

(ii) Assumption 2 in Theorem 1 implies that the functions $g_i, 0 \le i \le q$, are continuous and bounded on $]0, 1[$, and so they are Riemann integrable on this interval.

(iii) In the case $\gamma > 1$, the assumption (8) is a consequence of the assumption (7) and of Karamata's representation of positive slowly varying functions (see Bingham, Goldie and Teugels [7], Theorem 1.3.1).

**Proof of Theorem 1.** We set $\nu_0 = 0$ and fix the convention that $\phi(h)^{\nu_0} = 1$. Moreover, in the entire proof $K$ denotes a positive constant whose value does not matter. First we assume that $X$ is centered.

We prove the following asymptotic expansion for the expectation of $V_n(X)$:

$$\frac{n^{1-\gamma}}{L(1/n)} \mathbb{E}V_n(X) \overset{n \to +\infty}{=} \sum_{i=0}^q \left(\int_0^1 g_i(x)\,dx\right)\phi\left(\frac{1}{n}\right)^{\nu_i} + o\left(\phi\left(\frac{1}{n}\right)^{\nu_q}\right). \quad (10)$$

We have

$$d_{jk} = (\delta_1^{1/n} \circ \delta_2^{1/n} R)\left(\frac{j}{n}, \frac{k}{n}\right) \quad (11)$$

and

$$\mathbb{E}V_n(X) = \sum_{k=1}^{n-1} d_{kk}. \quad (12)$$

Moreover, the assumption (6) yields

$$\sup_{k=1,\ldots,n-1} \left| \frac{d_{kk}}{n^{\gamma-2}L(1/n)} - \sum_{i=0}^q g_i\left(\frac{k}{n}\right)\phi\left(\frac{1}{n}\right)^{\nu_i} \right| \overset{n \to +\infty}{=} o\left(\phi\left(\frac{1}{n}\right)^{\nu_q}\right). \quad (13)$$

Therefore,

$$\limsup_{n \to +\infty} \frac{1}{\phi(1/n)^{\nu_q}} \left| \frac{n^{1-\gamma}}{L(1/n)} \mathbb{E}V_n(X) - \sum_{i=0}^q \int_0^1 g_i(x)\,dx\,\phi\left(\frac{1}{n}\right)^{\nu_i} \right|$$

$$\leq \limsup_{n \to +\infty} \frac{1}{\phi(1/n)^{\nu_q}} \left| \sum_{k=1}^{n-1} \left( \frac{d_{kk}}{n^{\gamma-1}L(1/n)} - \frac{1}{n}\sum_{i=0}^{q} g_i\left(\frac{k}{n}\right)\phi\left(\frac{1}{n}\right)^{\nu_i} \right) \right|$$

$$+ \limsup_{n \to +\infty} \frac{1}{\phi(1/n)^{\nu_q}} \left| \sum_{i=0}^{q} \left( \int_0^1 g_i(x)\,\mathrm{d}x - \frac{1}{n}\sum_{k=1}^{n-1} g_i\left(\frac{k}{n}\right) \right)\phi\left(\frac{1}{n}\right)^{\nu_i} \right|$$

$$= L_1 + L_2.$$

We get

$$L_1 = \limsup_{n \to +\infty} \frac{1}{\phi(1/n)^{\nu_q}} \left| \sum_{k=1}^{n-1} \left( \frac{d_{kk}}{n^{\gamma-1}L(1/n)} - \frac{1}{n}\sum_{i=0}^{q} g_i\left(\frac{k}{n}\right)\phi\left(\frac{1}{n}\right)^{\nu_i} \right) \right|$$

$$\leq \limsup_{n \to +\infty} \frac{1}{\phi(1/n)^{\nu_q}} \sup_{k=1,\dots,n-1} \left| \frac{d_{kk}}{n^{\gamma-2}L(1/n)} - \sum_{i=0}^{q} g_i\left(\frac{k}{n}\right)\phi\left(\frac{1}{n}\right)^{\nu_i} \right|.$$

Thus (13) implies that $L_1 = 0$.

For $L_2$, notice that

$$L_2 = \limsup_{n \to +\infty} \frac{1}{\phi(1/n)^{\nu_q}} \left| \sum_{i=0}^{q} \left( \int_0^1 g_i(x)\,\mathrm{d}x - \frac{1}{n}\sum_{k=1}^{n-1} g_i\left(\frac{k}{n}\right) \right)\phi\left(\frac{1}{n}\right)^{\nu_i} \right|$$

$$\leq \limsup_{n \to +\infty} \frac{1}{\phi(1/n)^{\nu_q}} \sum_{i=0}^{q} \left| \int_0^{(n-1)/n} g_i(x)\,\mathrm{d}x - \frac{1}{n}\sum_{k=1}^{n-1} g_i\left(\frac{k}{n}\right) \right| \phi\left(\frac{1}{n}\right)^{\nu_i}$$

$$+ \sum_{i=0}^{q} \limsup_{n \to +\infty} \frac{\phi(1/n)^{\nu_i}}{\phi(1/n)^{\nu_q}} \int_{(n-1)/n}^{t} g_i(x)\,\mathrm{d}x = L_2^{(1)} + L_2^{(2)}.$$

The term $L_2^{(2)}$ is obviously equal to 0 due to (7) and the fact that the functions $g_i$ are bounded.

Moreover the assumption 2(c)(i) in Theorem 1 implies that there exists $K > 0$ such that for all $0 \leq i \leq q-1$,

$$\left| \int_0^{(n-1)/n} g_i(x)\,\mathrm{d}x - \frac{1}{n}\sum_{k=1}^{n-1} g_i\left(\frac{k}{n}\right) \right|$$

$$\leq \sum_{k=1}^{n-1} \int_{(k-1)/n}^{k/n} \left| g_i(x) - g_i\left(\frac{k}{n}\right) \right| \mathrm{d}x \leq K \sum_{k=1}^{n-1} \int_{(k-1)/n}^{k/n} \left( \frac{k}{n} - x \right) \mathrm{d}x \leq \frac{K}{n}.$$

Consequently,

$$L_2^{(1)} = \limsup_{n \to +\infty} \frac{1}{\phi(1/n)^{\nu_q}} \sum_{i=0}^{q} \left| \int_0^{(n-1)/n} g_i(x)\,\mathrm{d}x - \frac{1}{n}\sum_{k=1}^{n-1} g_i\left(\frac{k}{n}\right) \right| \phi\left(\frac{1}{n}\right)^{\nu_i}$$

$$\leq \limsup_{n \to +\infty} \frac{1}{\phi(1/n)^{\nu_q}} \sum_{i=0}^{q-1} \left| \int_0^{(n-1)/n} g_i(x)\, dx - \frac{1}{n} \sum_{k=1}^{n-1} g_i\left(\frac{k}{n}\right) \right| \phi\left(\frac{1}{n}\right)^{\nu_i}$$

$$+ \limsup_{n \to +\infty} \left| \int_0^{(n-1)/n} g_q(x)\, dx - \frac{1}{n} \sum_{k=1}^{n-1} g_q\left(\frac{k}{n}\right) \right|$$

$$\leq K \sum_{i=0}^{q-1} \limsup_{n \to +\infty} \frac{1}{n\phi(1/n)^{\nu_q}} \phi\left(\frac{1}{n}\right)^{\nu_i}$$

$$+ \limsup_{n \to +\infty} \left| \int_0^{(n-1)/n} g_q(x)\, dx - \frac{1}{n} \sum_{k=1}^{n-1} g_q\left(\frac{k}{n}\right) \right|,$$

where the first term of the right-hand side is equal to 0 because of (7) and the assumption 2(c)(i) in Theorem 1, and the second term is equal to 0 according to classical results on Riemann sums. Therefore, $L_2^{(1)} = 0$. This proves the asymptotic expansion (10).

Next we prove that almost surely

$$\frac{n^{1-\gamma}}{L(1/n)}(V_n(X) - \mathbb{E}V_n(X)) \overset{n \to \pm \infty}{=} o\left( \phi\left(\frac{1}{n}\right)^{\nu_q} \right). \tag{14}$$

Application of Cochran's theorem to the Gaussian vector yields

$$\sqrt{\frac{n^{1-\gamma}}{\phi(1/n)^{\nu_q} L(1/n)}} \Delta X_k,$$

So there are $n-1$ nonnegative real numbers $(\mu_{1,n}, \ldots, \mu_{n-1,n})$ and one $(n-1)$-dimensional Gaussian vector $Y_n$, such that its components are independent Gaussian variables $\mathcal{N}(0,1)$ and

$$\frac{n^{1-\gamma}}{L(1/n)} V_n(X) = \sum_{j=1}^{n-1} \mu_{j,n}(Y_n^{(j)})^2. \tag{15}$$

As in Bégyn [5], Hanson and Wright's inequality (see Hanson and Wright [12]) yields that, for all $0 < \varepsilon < 1$,

$$\mathbb{P}\left( \frac{n^{1-\gamma}}{L(1/n)\phi(1/n)^{\nu_q}} |V_n(X) - \mathbb{E}V_n(X)| \geq \varepsilon \right) \leq 2\exp\left( -K\varepsilon^2 n\phi\left(\frac{1}{n}\right)^{\nu_q} \right). \tag{16}$$

So if we set

$$\varepsilon_n^2 = \frac{2\log n}{Kn} \phi\left(\frac{1}{n}\right)^{-\nu_q},$$

it follows from (7) that

$$\lim_{n \to +\infty} \varepsilon_n = 0 \quad \text{and} \quad \sum_{n=0}^{+\infty} \mathbb{P}\left(\frac{n^{1-\gamma}}{L(1/n)\phi(1/n)^{\nu_q}}|V_n(X) - \mathbb{E}V_n(X)| \geq \varepsilon_n\right) < +\infty,$$

and the Borel–Cantelli lemma yields (14).

Now let us examine the case of non-centered $X$. Set $M = \{M_t; t \in [0,1]\}$. From assumption 1,

$$\frac{n^{1-\gamma}}{L(1/n)}V_n(M) \stackrel{n \to +\infty}{=} \mathcal{O}\left(\frac{1}{n^\gamma L(1/n)}\right),$$

and by adding (8) we obtain

$$\lim_{n \to +\infty} \frac{n^{1-\gamma}}{L(1/n)}V_n(M) = 0.$$

So, if we apply the theorem to the centered process $\widetilde{X}_t = X_t - \mathbb{E}(X_t)$, using the arguments of Baxter [3], we obtain the result for $X$. $\qquad\qquad\square$

In the sequel, we apply these results to the identification of some fractional models. We will obtain strongly consistent estimators that will be more interesting in practice if they are asymptotically normal. Therefore, we establish a central limit theorem for $V_n(X)$.

## 2.2. Central limit theorem

The integral

$$\rho_\gamma(j,k) = \int_j^{j+1} \mathrm{d}u \int_{u-1}^u \mathrm{d}v \int_k^{k+1} \mathrm{d}x \int_{x-1}^x \frac{1}{(v-y)^{2+\gamma}}\,\mathrm{d}y \tag{17}$$

with $j - k \geq 2$ is absolutely convergent when $\gamma < 2$. Because it depends only on the difference $j - k$, we denote it $\rho_\gamma(j-k)$.

By considering $l = j - k \geq 2$ and $0 < \gamma < 2$, we obtain the following equalities: If $\gamma \neq 1$,

$$\rho_\gamma(l) = \frac{(|l-2|^{2-\gamma} - 4|l-1|^{2-\gamma} + 6|l|^{2-\gamma} - 4|l+1|^{2-\gamma} + |l+2|^{2-\gamma})}{(\gamma-2)(\gamma-1)\gamma(\gamma+1)}; \tag{18}$$

if $\gamma = 1$,

$$\begin{aligned}
\rho_1(l) = \tfrac{1}{2}(&|l-2|\log|l-2| - 4|l-1|\log|l-1| + 6|l|\log|l| \\
&- 4|l+1|\log|l+1| + |l+2|\log|l+2|).
\end{aligned} \tag{19}$$

Moreover, we notice that (17) yields the existence of a constant $K > 0$ such that, for all $l \geq 2$, we have $|\rho_\gamma(l)| \leq K l^{-2-\gamma}$. For $\gamma \in\, ]0, 2[$, set

$$\|\rho_\gamma\|^2 = \sum_{l=2}^{+\infty} \rho_\gamma(l)^2. \tag{20}$$

We may now prove a central limit theorem with additional assumptions. The preceding formulas will be useful to compute the asymptotic behavior of $d_{jk}$.

**Theorem 2.** *Assume that $X$ is centered and satisfies the following statements:*

1. *$R$ is continuous in $[0, 1]^2$.*
2. *Let $T = \{0 \leq t \leq s \leq 1\}$. We assume that the derivative $\frac{\partial^4 R}{\partial s^2 \partial t^2}$ exists in $]0, 1]^2 \setminus \{s = t\}$, and that there exist a continuous function $C : T \to \mathbb{R}$, a real $\gamma \in\, ]0, 2[$ and a positive slowly varying function $L : ]0, 1] \to \mathbb{R}$ such that*

$$\forall (s, t) \in \overset{\circ}{T} \qquad \frac{(s-t)^{2+\gamma}}{L(s-t)} \frac{\partial^4 R}{\partial s^2 \partial t^2}(s, t) = C(s, t), \tag{21}$$

   *where $\overset{\circ}{T}$ denotes the interior of $T$ (i.e., $\overset{\circ}{T} = \{0 < t < s < 1\}$).*
3. *We assume that there exist $q + 1$ functions ($q \in \mathbb{N}$) $g_0, g_1, \ldots, g_q$ from $]0, 1[$ to $\mathbb{R}$, $q$ real numbers $0 < \nu_1 < \cdots < \nu_q$ and a function $\phi : ]0, 1[ \to ]0, +\infty[$ such that:*
   (a) *if $q \geq 1$, then $\forall 0 \leq i \leq q-1$, $g_i$ is Lipschitz on $]0, 1[$;*
   (b) *$g_q$ is $(1/2 + \alpha_q)$-Hölderian on $]0, 1[$ with $0 < \alpha_q \leq 1/2$;*
   (c) *there exists $t \in\, ]0, 1[$ such that $g_0(t) \neq 0$;*
   (d) *we have*

$$\lim_{h \to 0^+} \frac{1}{\sqrt{h}} \left( \sup_{h \leq t \leq 1-h} \left| \frac{(\delta_1^h \circ \delta_2^h R)(t, t)}{h^{2-\gamma} L(h)} - g_0(t) - \sum_{i=1}^q g_i(t) \phi(h)^{\nu_i} \right| \right) = 0, \tag{22}$$

   *where if $q = 0$, then $\sum_{i=1}^q g_i(t) \phi(h)^{\nu_i} = 0$, and where if $q \neq 0$, then $\lim_{h \to 0^+} \phi(h) = 0$;*
   (e) *there exists a bounded function $\widetilde{g} : ]0, 1[ \to \mathbb{R}$ such that*

$$\lim_{h \to 0^+} \sup_{h \leq t \leq 1-2h} \left| \frac{(\delta_1^h \circ \delta_2^h R)(t+h, t)}{h^{2-\gamma} L(h)} - \widetilde{g}(t) \right| = 0. \tag{23}$$

*Then we have*

$$\sqrt{n} \left( \frac{n^{1-\gamma}}{L(1/n)} V_n(X) - \int_0^1 g_0(x)\,dx - \sum_{i=1}^q \int_0^1 g_i(x)\,dx \cdot \phi\left(\frac{1}{n}\right)^{\nu_i} \right) \xrightarrow{(\mathcal{L})} \mathcal{N}(0, \sigma^2), \tag{24}$$

*where*

$$\sigma^2 = 2 \int_0^1 g_0(x)^2\,dx + 4 \int_0^1 \widetilde{g}(x)^2\,dx + 4 \|\rho_\gamma\|^2 \int_0^1 C(x, x)^2\,dx. \tag{25}$$

**Remark.** (i) In Theorem 1, there was no minimum value for the integer $q$, but this is not the case in the assumption (22): we must choose $q$ large enough such that the bias is negligible with respect to the stochastic error.

(ii) Assumption (22) yields that the functions $g_i$, $0 \leq i \leq q$, are continuous on $]0, 1[$.

**Proof of Theorem 2.** In all the proof, $K$ denotes a positive constant whose value does not matter. To simplify notation, choose the convention $\nu_0 = 0$ and $\forall h \in ]0, 1[$, $\phi(h)^{\nu_0} = 1$. Set

$$b_n = \sum_{i=0}^{q} \int_0^1 g_i(x) \, dx \, \phi\left(\frac{1}{n}\right)^{\nu_i}, \qquad T_n = \sqrt{n} \, \frac{n^{1-\gamma}}{L(1/n)} V_n(X), \qquad \widetilde{T_n} = T_n - \mathbb{E}(T_n). \quad (26)$$

We split the proof into three steps: in the first and second steps, we prove the convergence when $n \to +\infty$ of $\widetilde{T_n}$ toward a centered Gaussian law with variance $\sigma^2$; in the third step, we prove the conclusion of Theorem 2.

*Step 1.* We note that $V_n(X)$ is the square of the Euclidean norm of a $(n-1)$-dimensional Gaussian vector whose components are

$$\sqrt{\frac{n^{1-\gamma}}{L(1/n)}} \Delta X_k, \qquad 1 \leq k \leq n-1.$$

Hence, by the classical Cochran theorem, we can find $a_n = n-1$ positive real numbers $(\lambda_{1,n}, \ldots, \lambda_{a_n,n})$ and one $a_n$-dimensional Gaussian vector $Y_n$, such that its components are independent Gaussian variables $\mathcal{N}(0, 1)$ and

$$V_n(X) = \sum_{j=1}^{a_n} \lambda_{j,n} (Y_n^{(j)})^2, \quad (27)$$

with the convention that the empty sum is equal to zero.

We set

$$S_n(X) = V_n(X) - \mathbb{E} V_n(X).$$

We want to apply the Lindeberg central limit theorem to $S_n(X)$. We must verify that

$$\lambda_n^* = \max_{1 \leq j \leq a_n} \lambda_{j,n} = o(\sqrt{\operatorname{Var} S_n(X)}). \quad (28)$$

We have

$$\lambda_n^* \leq K \frac{n^{1-\gamma}}{L(1/n)} \max_{1 \leq k \leq n-1} \sum_{j=1}^{n-1} |d_{jk}|.$$

With the same methods as in Bégyn [5], we can check that

$$\lambda_n^* \overset{n \to +\infty}{=} \mathcal{O}\left(\frac{1}{n}\right).$$

We have

$$\mathbb{E}[V_n(X)^2] = \sum_{j,k=1}^{n-1} \mathbb{E}[(\Delta X_j)^2 (\Delta X_k)^2].$$

Because the vector $(\Delta X_k)_{1 \leq k \leq N-1}$ is Gaussian, Isserlis formulas yield (see Isserlis [14])

$$\mathbb{E}[(\Delta X_k)^4] = 3(\mathbb{E}[(\Delta X_k)^2])^2$$

and, if $j \neq k$,

$$\mathbb{E}[(\Delta X_j)^2 (\Delta X_k)^2] = \mathbb{E}[(\Delta X_j)^2]\mathbb{E}[(\Delta X_k)^2] + 2(\mathbb{E}[\Delta X_j \Delta X_k])^2.$$

Therefore,

$$\operatorname{Var} V_n(X) = 2 \sum_{k=1}^{n-1} d_{kk}^2 + 4 \sum_{1 \leq k < j \leq n-1} d_{jk}^2 \tag{29}$$

and, consequently,

$$\operatorname{Var} S_n(X) \geq 2 \sum_{k=1}^{n-1} d_{kk}^2.$$

Moreover, the assumption (22) yields

$$\lim_{n \to +\infty} \sup_{1 \leq k \leq n-1} \left| \frac{(\delta_1^{1/n} \circ \delta_2^{1/n} R)(k/n, k/n)}{n^{\gamma-2} L(1/n)} - g_0\left(\frac{k}{n}\right) \right| = 0$$

and, because $g_0$ is bounded on $]0,1[$,

$$\lim_{n \to +\infty} \sup_{1 \leq k \leq n-1} \left| \left( \frac{(\delta_1^{1/n} \circ \delta_2^{1/n} R)(k/n, k/n)}{n^{\gamma-2} L(1/n)} \right)^2 - g_0\left(\frac{k}{n}\right)^2 \right| = 0.$$

With the same ideas as in the proof of (10), we can show that the previous limit yields

$$\lim_{n \to +\infty} \frac{n^{3-2\gamma}}{L(1/n)^2} \sum_{k=1}^{n-1} d_{kk}^2 = \int_0^1 g_0^2(x) \, \mathrm{d}x.$$

Therefore,

$$\liminf_{n \to +\infty} \frac{n^{3-2\gamma}}{L(1/n)^2} \operatorname{Var} S_n(X) \geq 2 \int_0^1 g_0^2(x) \, \mathrm{d}x > 0.$$

Thus there exists $K > 0$ such that

$$0 \leq \frac{\lambda_n^*}{\sqrt{\operatorname{Var} S_n(X)}} \leq \frac{K}{\sqrt{n}},$$

which yields (28).

Consequently, the Lindeberg central limit theorem yields that when $n \to +\infty$,

$$\frac{S_n(X)}{\sqrt{\operatorname{Var} S_n(X)}} \xrightarrow{(\mathcal{L})} \mathcal{N}(0, 1).$$

Equivalently, we have shown that when $n \to +\infty$,

$$\frac{\widetilde{T_n}}{\sqrt{\operatorname{Var} \widetilde{T_n}}} \xrightarrow{(\mathcal{L})} \mathcal{N}(0, 1). \tag{30}$$

*Step 2.* Let us prove that

$$\lim_{n \to +\infty} \operatorname{Var} \widetilde{T_n} = \sigma, \tag{31}$$

where $\sigma^2$ was defined in (25).

We have seen (29):

$$\operatorname{Var} S_n(X) = 2 \sum_{k=1}^{n-1} d_{kk}^2 + 4 \sum_{1 \le k < j \le n-1} d_{jk}^2.$$

*Step 2.1.* Let us prove that

$$\lim_{n \to +\infty} \frac{n^{3-2\gamma}}{L(1/n)^2} \sum_{\substack{1 \le k \le j \le n-1 \\ j-k \ge 3}} d_{jk}^2 = \int_0^1 C(x,x)^2 \, \mathrm{d}x \sum_{l=3}^{+\infty} \rho_\gamma(l)^2. \tag{32}$$

If $j \neq 1$, $k \neq 1$ and $j - k \ge 3$, then

$$d_{jk} = (\delta_1^{1/n} \circ \delta_2^{1/n}) R\left(\frac{j}{n}, \frac{k}{n}\right)$$

$$= \int_{j/n}^{(j+1)/n} \mathrm{d}u \int_{u-(1/n)}^u \mathrm{d}v \int_{k/n}^{(k+1)/n} \mathrm{d}x \int_{x-(1/n)}^x \frac{\partial^4 R}{\partial s^2 \, \partial t^2}(v, y) \, \mathrm{d}y$$

$$= \int_{j/n}^{(j+1)/n} \mathrm{d}u \int_{u-(1/n)}^u \mathrm{d}v \int_{k/n}^{(k+1)/n} \mathrm{d}x \int_{x-(1/n)}^x \frac{C(v,y)}{(v-y)^{2+\gamma}} L(v-y) \, \mathrm{d}y.$$

We set

$$\varepsilon_n = \sup\left\{ \left| C(v,y) - C\left(\frac{j}{n}, \frac{k}{n}\right) \right|; 3 \le j - k \le n-2, \frac{j-1}{n} \le v \le \frac{j+1}{n}, \frac{k-1}{n} \le y \le \frac{k+1}{n} \right\}.$$

Because $C$ is uniformly continuous on the compact set $T$,

$$\lim_{n \to +\infty} \varepsilon_n = 0.$$

Moreover, we set

$$r_n(j,k) = \int_{j/n}^{(j+1)/n} \mathrm{d}u \int_{u-(1/n)}^{u} \mathrm{d}v \int_{k/n}^{(k+1)/n} \mathrm{d}x \int_{x-(1/n)}^{x} \frac{L(v-y)}{(v-y)^{2+\gamma}} \, \mathrm{d}y$$

$$= n^{\gamma-2} \int_{j}^{j+1} \mathrm{d}u \int_{u-1}^{u} \mathrm{d}v \int_{k}^{k+1} \mathrm{d}x \int_{x-1}^{x} \frac{L((v-y)/n)}{(v-y)^{2+\gamma}} \, \mathrm{d}y,$$

and because it depends only on $j-k$, we denote it $r_n(j-k)$.

We have

$$\left| d_{jk} - C\left(\frac{j}{n}, \frac{k}{n}\right) r_n(j-k) \right|$$

$$= \left| \int_{j/n}^{(j+1)/n} \mathrm{d}u \int_{u-(1/n)}^{u} \mathrm{d}v \int_{k/n}^{(k+1)/n} \mathrm{d}x \int_{x-(1/n)}^{x} \frac{C(v,y) - C(j/n, k/n)}{(v-y)^{2+\gamma}} L(v-y) \, \mathrm{d}y \right|$$

$$\leq \varepsilon_n r_n(j-k).$$

So we must find an upper bound for $r_n(j-k)$. Let us note that the function $\psi(h) = h^{-\gamma/2} L(h)$ is regularly varying of index $-\gamma/2 < 0$. Therefore, the Karamata theorem of uniform convergence (see Bingham *et al.* [7], Theorem 1.5.2) yields

$$\lim_{h \to 0^+} \frac{\psi(zh)}{\psi(h)} = \frac{1}{z^{\gamma/2}}$$

uniformly in $z \in [3, +\infty[$.

As a consequence, there exists $K > 0$ (which depends only on $L$ and $\gamma$) such that for $n$ large enough,

$$\forall z \geq 3, \qquad \psi\left(\frac{z}{n}\right) \leq K\psi\left(\frac{1}{n}\right). \tag{33}$$

However, we have

$$r_n(j-k) = n^{(\gamma/2)-2} \int_{j}^{j+1} \mathrm{d}u \int_{u-1}^{u} \mathrm{d}v \int_{k}^{k+1} \mathrm{d}x \int_{x-1}^{x} \frac{\psi((v-y)/n)}{(v-y)^{2+(\gamma/2)}} \, \mathrm{d}y.$$

Therefore, (17) yields

$$r_n(j-k) \leq Kn^{\gamma-2} L\left(\frac{1}{n}\right) \rho_{\gamma/2}(j-k) \leq K \frac{n^{\gamma-2} L(1/n)}{(j-k-2)^{2+(\gamma/2)}}.$$

Consequently, we have

$$\left| d_{jk} - C\left(\frac{j}{n}, \frac{k}{n}\right) r_n(j-k) \right| \leq K\varepsilon_n \frac{n^{\gamma-2} L(1/n)}{(j-k-2)^{2+(\gamma/2)}},$$

and because $C$ is bounded,

$$\left| d_{jk}^2 - C\left(\frac{j}{n}, \frac{k}{n}\right)^2 r_n(j-k)^2 \right| \le K\varepsilon_n \frac{n^{2\gamma-4}L(1/n)^2}{(j-k-2)^{4+\gamma}}.$$

Using the same perturbation argument as in Bégyn ([5], pages 10–11), we can check that it is still true whenever $j=1$ or $k=1$. Consequently,

$$\frac{n^{3-2\gamma}}{L(1/n)^2}\left| \sum_{l=3}^{n-2}\sum_{k=1}^{n-1-l} d_{l+k,k}^2 - \sum_{l=3}^{n-2} r_n(l)^2 \sum_{k=1}^{n-1-l} C\left(\frac{k+l}{n}, \frac{k}{n}\right)^2 \right|$$

$$\le K\varepsilon_n \frac{1}{n}\sum_{l=3}^{n-2}\sum_{k=1}^{n-1-l} \frac{1}{(l-2)^{4+\gamma}} \le K\varepsilon_n \sum_{l=3}^{+\infty} \frac{1}{(l-2)^{4+\gamma}}$$

and the right-hand side is convergent because $4+\gamma > 1$. This yields

$$\lim_{n\to+\infty} \frac{n^{3-2\gamma}}{L(1/n)^2}\left( \sum_{l=3}^{n-2}\sum_{k=1}^{n-1-l} d_{l+k,k}^2 - \sum_{l=3}^{n-2} r_n(l)^2 \sum_{k=1}^{n-1-l} C\left(\frac{k+l}{n}, \frac{k}{n}\right)^2 \right) = 0. \qquad (34)$$

Moreover,

$$\lim_{n\to+\infty} \frac{1}{n} \sum_{k=1}^{n-1-l} C\left(\frac{k+l}{n}, \frac{k}{n}\right)^2 = \lim_{n\to+\infty} \frac{1}{n} \sum_{k=1}^{n-1-l} C\left(\frac{k}{n}, \frac{k}{n}\right)^2 = \int_0^1 C(x,x)^2\,\mathrm{d}x,$$

thanks to the uniform continuity of $C$.

In addition, because $L$ is slowly varying, the theorem of dominated convergence yields (using inequality (33))

$$\lim_{n\to+\infty} \frac{r_n(l)}{n^{\gamma-2}L(1/n)} = \rho_\gamma(l).$$

Hence for all $l \ge 3$,

$$\lim_{n\to+\infty} \frac{r_n(l)^2}{n^{2\gamma-4}L(1/n)}\frac{1}{n} \sum_{k=1}^{n-1-l} C\left(\frac{k+l}{n}, \frac{k}{n}\right)^2 = \rho_\gamma(l)^2 \int_0^1 C(x,x)^2\,\mathrm{d}x,$$

and (using (33)) we can check that $\forall l \ge 3$,

$$\frac{r_n(l)^2}{n^{2\gamma-4}L(1/n)^2}\frac{1}{n} \sum_{k=1}^{n-1-l} C\left(\frac{k+l}{n}, \frac{k}{n}\right)^2 \le K\rho_{\gamma/2}(l)^2 \le \frac{K}{(l-2)^{4+\gamma}}.$$

Consequently, the theorem of dominated convergence for series yields

$$\lim_{n\to+\infty} \frac{n^{3-2\gamma}}{L(1/n)^2} \sum_{l=3}^{n-2} r_n(l)^2 \sum_{k=1}^{n-1-l} C\left(\frac{k+l}{n}, \frac{k}{n}\right)^2 = \int_0^1 C(x,x)^2\,\mathrm{d}x \sum_{l=3}^{+\infty} \rho_\gamma(l)^2.$$

With (34), we obtain

$$\lim_{n \to +\infty} \frac{n^{3-2\gamma}}{L(1/n)^2} \sum_{\substack{1 \le k < j \le n-1 \\ j-k \ge 3}} d_{jk}^2 = \lim_{n \to +\infty} \frac{n^{3-2\gamma}}{L(1/n)^2} \sum_{l=3}^{n-2} \sum_{k=1}^{n-1-l} d_{l+k,k}^2$$

$$= \int_0^1 C(x,x)^2 \, \mathrm{d}x \sum_{l=3}^{+\infty} \rho_\gamma(l)^2.$$

*Step 2.2.* Let us prove that

$$\lim_{n \to +\infty} \frac{n^{3-2\gamma}}{L(1/n)^2} \sum_{\substack{1 \le k < j \le n-1 \\ j-k=2}} d_{jk}^2 = \rho_\gamma(2)^2 \int_0^1 C(x,x)^2 \, \mathrm{d}x. \qquad (35)$$

With the perturbation argument of Bégyn [5], we can check that

$$\left| d_{k+2,k} - C\left(\frac{k+2}{n}, \frac{k}{n}\right) r_n(2) \right| \le K \varepsilon_n n^{\gamma-2} L\left(\frac{1}{n}\right) \rho_{\gamma/2}(2)$$

and, consequently,

$$\left| d_{k+2,k}^2 - C\left(\frac{k+2}{n}, \frac{k}{n}\right)^2 r_n(2)^2 \right| \le K \varepsilon_n n^{2\gamma-4} L\left(\frac{1}{n}\right)^2 \rho_{\gamma/2}(2)^2.$$

Therefore, using the same arguments as in Step 2.1, we have

$$\lim_{n \to +\infty} \frac{n^{3-2\gamma}}{L(1/n)^2} \sum_{\substack{1 \le k < j \le n-1 \\ j-k=2}} d_{jk}^2 = \lim_{n \to +\infty} \frac{n^{3-2\gamma}}{L(1/n)^2} \sum_{k=1}^{n-3} d_{k+2,k}^2$$

$$= \lim_{n \to +\infty} \frac{r_n(2)^2}{n^{2\gamma-4} L(1/n)} \frac{1}{n} \sum_{k=1}^{n-3} C\left(\frac{k+2}{n}, \frac{k}{n}\right)^2$$

$$= \rho_\gamma(2)^2 \int_0^1 C(x,x)^2 \, \mathrm{d}x.$$

*Step 2.3.* Let us prove that

$$\lim_{n \to +\infty} \frac{n^{3-2\gamma}}{L(1/n)^2} \sum_{\substack{1 \le k < j \le n-1 \\ j-k=1}} d_{jk}^2 = \int_0^1 \widetilde{g}(x)^2 \, \mathrm{d}x. \qquad (36)$$

This is a consequence of the assumption (23).

*Step 2.4.* Let us prove that

$$\lim_{n \to +\infty} \frac{n^{3-2\gamma}}{L(1/n)^2} \sum_{k=1}^{n-1} d_{kk}^2 = \int_s^t g_0(x)^2 \, \mathrm{d}x. \tag{37}$$

This is a consequence of

$$\lim_{h \to 0^+} \sup_{h \le t \le 1-h} \left| \frac{(\delta_1^h \circ \delta_2^h R)(t,t)}{L(h)h^{2-\gamma}} - g_0(t) \right| = 0,$$

which comes from the assumption (22).

The preceding four steps imply (31). Let us remark that (30) and (31) yield that

$$\widetilde{T_n} \xrightarrow{(\mathcal{L})} \mathcal{N}(0, \sigma). \tag{38}$$

*Step 3.* To prove Theorem 2, we use the decomposition

$$T_n - \sqrt{n} b_n = \widetilde{T_n} + \mathbb{E}T_n - \sqrt{n} b_n. \tag{39}$$

Let us prove that

$$\lim_{n \to +\infty} |\mathbb{E}T_n - \sqrt{n}\, b_n| = 0. \tag{40}$$

We have

$$\mathbb{E}(T_n) - \sqrt{n}\, b_n = \sqrt{n} \left( \frac{n^{1-\gamma}}{L(1/n)} \sum_{k=1}^{n-1} d_{kk} - \sum_{i=0}^{q} \int_0^t g_i(x) \, \mathrm{d}x \, \phi \left( \frac{1}{n} \right)^{\nu_i} \right).$$

On the one hand,

$$\sqrt{n} \left| \frac{n^{1-\gamma}}{L(1/n)} \sum_{k=1}^{n-1} d_{kk} - \frac{1}{n} \sum_{i=0}^{q} \phi \left( \frac{1}{n} \right)^{\nu_i} \sum_{k=1}^{n-1} g_i \left( \frac{k}{n} \right) \right|$$

$$\le \sqrt{n} \sup_{k=1,\dots,n-1} \left| \frac{d_{kk}}{n^{\gamma-2}L(1/n)} - \sum_{i=0}^{q} g_i \left( \frac{k}{n} \right) \phi \left( \frac{1}{n} \right)^{\nu_i} \right|.$$

Therefore, the assumption (22) yields

$$\lim_{n \to +\infty} \sqrt{n} \left| \frac{n^{1-\gamma}}{L(1/n)} \sum_{k=1}^{n-1} d_{kk} - \frac{1}{n} \sum_{i=0}^{q} \phi \left( \frac{1}{n} \right)^{\nu_i} \sum_{k=1}^{n-1} g_i \left( \frac{k}{n} \right) \right| = 0. \tag{41}$$

Moreover, if we choose $0 \le i \le q-1$, we have

$$\sqrt{n} \left| \frac{1}{n} \phi \left( \frac{1}{n} \right)^{\nu_i} \sum_{k=1}^{n-1} g_i \left( \frac{k}{n} \right) - \phi \left( \frac{1}{n} \right)^{\nu_i} \int_0^1 g_i(x) \, \mathrm{d}x \right| \le K \phi \left( \frac{1}{n} \right)^{\nu_i} \frac{1}{\sqrt{n}},$$

due to the fact that $g_i$ is bounded and Lipschitz on $]0, 1[$ (see the Proof of Theorem 1).

Furthermore,

$$\sqrt{n}\left|\frac{1}{n}\phi\left(\frac{1}{n}\right)^{\nu_q}\sum_{k=1}^{n-1}g_q\left(\frac{k}{n}\right) - \phi\left(\frac{1}{n}\right)^{\nu_q}\int_0^1 g_q(x)\,\mathrm{d}x\right| \leq K\phi\left(\frac{1}{n}\right)^{\nu_q}\left(\frac{1}{n}\right)^{\nu_q},$$

knowing that $g_q$ is bounded and $1/2 + \alpha_q$-Hölderian on $]0, 1[$. Consequently,

$$\lim_{n\to+\infty}\sqrt{n}\left|\frac{1}{n}\sum_{i=0}^q\left(\phi\left(\frac{1}{n}\right)^{\nu_i}\sum_{k=1}^{n-1}g_i\left(\frac{k}{n}\right) - \int_0^1 g_i(x)\,\mathrm{d}x\,\phi\left(\frac{1}{n}\right)^{\nu_i}\right)\right| = 0. \tag{42}$$

To finish the proof, just note that (40) is a consequence of (41) and (42).

Next, by combining the Prokhorov theorem, all the preceding steps and the Slutzky lemma with (39) and (40), we get (24). □

In the sequel, we consider estimators of some functions of the parameters $\gamma$, which are constructed with both second-order quadratic variations $V_n(X)$ and $V_{2n}(X)$. So the preceding theorem is not sufficient to prove the asymptotic normality of the estimators.

## 2.3. Bivariate central limit theorem

The next theorem will be useful to prove the asymptotic normality of our estimators. We define the following constants, which appear in the asymptotic covariance of the two quadratic variations:

$$\sigma_{1,\mathrm{cov}}^2 = 2\int_0^1 \widetilde{g}(t)^2\,\mathrm{d}t + 4\rho_\gamma(2)\int_0^1 g_0(t)C(t,t)\,\mathrm{d}t + 4\rho_\gamma(3)\int_0^1 \widetilde{g}(t)C(t,t)\,\mathrm{d}t$$

$$+ 4\int_0^1 C(t,t)^2\,\mathrm{d}t\sum_{l=4}^{+\infty}\rho_\gamma(l)\rho_\gamma(l-2),$$

$$\sigma_{2,\mathrm{cov}}^2 = 4\int_0^1 g_0(t)\widetilde{g}(t)\,\mathrm{d}t + 4\rho_\gamma(2)\int_0^1 \widetilde{g}(t)C(t,t)\,\mathrm{d}t$$

$$+ 4\int_0^1 C(t,t)^2\,\mathrm{d}t\sum_{l=3}^{+\infty}\rho_\gamma(l)\rho_\gamma(l-1)$$

and

$$\sigma_*^2 = 3\sigma^2 + \sigma_{1,\mathrm{cov}}^2 + 4\sigma_{2,\mathrm{cov}}^2. \tag{43}$$

**Theorem 3.** *We use the same assumptions as in Theorem [2](#). Then we have*

$$\sqrt{n} \left( \begin{array}{c} \frac{n^{1-\gamma}}{L(\frac{1}{n})} V_n(X) - \int_0^1 g_0(x) \, dx - \sum_{i=1}^q \int_0^1 g_i(x) \, dx \cdot \phi(\frac{1}{n})^{\nu_i} \\ \frac{(2n)^{1-\gamma}}{L(\frac{1}{2n})} V_{2n}(X) - \int_0^1 g_0(x) \, dx - \sum_{i=1}^q \int_0^1 g_i(x) \, dx \cdot \phi(\frac{1}{2n})^{\nu_i} \end{array} \right) \xrightarrow{(\mathcal{L})} \mathcal{N}(0, \Sigma), \quad (44)$$

*where the $2 \times 2$ matrix $\Sigma$ is defined by*

$$\Sigma = \begin{bmatrix} \sigma^2 & 2^{\gamma-2}\sigma_*^2 \\ 2^{\gamma-2}\sigma_*^2 & \sigma^2/2 \end{bmatrix}. \quad (45)$$

**Proof.** We set for $\lambda, \mu \in \mathbb{R}$,

$$S_n(\lambda, \mu) = \lambda \frac{n^{1-\gamma}}{L(1/n)} V_n(X) + \mu \frac{(2n)^{1-\gamma}}{L(1/(2n))} V_{2n}(X).$$

We begin by showing that when $\lambda, \mu \geq 0$,

$$\sqrt{n} \left( S_n(\lambda, \mu) - \lambda \sum_{i=0}^q \int_0^1 g_i(x) \, dx \, \phi\left(\frac{1}{n}\right)^{\nu_i} - \mu \sum_{i=0}^q \int_0^1 g_i(x) \, dx \cdot \phi\left(\frac{1}{2n}\right)^{\nu_i} \right)$$

$$\xrightarrow{(\mathcal{L})} \mathcal{N}\left( 0, \lambda^2 \sigma^2 + \mu^2 \frac{\sigma^2}{2} + 2^{\gamma-1} \lambda \mu \sigma_*^2 \right). \quad (46)$$

First we must prove that

$$\lim_{n \to +\infty} n \operatorname{Var} S_n(\lambda, \mu) = \lambda^2 \sigma^2 + \mu^2 \frac{\sigma^2}{2} + 2^{\gamma-1} \lambda \mu \sigma_*^2. \quad (47)$$

We have

$$\operatorname{Var} S_n(\lambda, \mu) = \lambda^2 \frac{n^{2-2\gamma}}{L(1/n)^2} \operatorname{Var} V_n(X) + \mu^2 \frac{(2n)^{2-2\gamma}}{L(1/(2n))^2} \operatorname{Var} V_{2n}(X)$$

$$+ 2\lambda\mu \frac{2^{1-\gamma} n^{2-2\gamma}}{L(1/n) L(1/(2n))} \operatorname{Cov}(V_n(X), V_{2n}(X)).$$

In the Proof of Theorem [2](#) we showed [(31)](#):

$$\lim_{n \to +\infty} \frac{n^{3-2\gamma}}{L(1/n)^2} \operatorname{Var} V_n(X) = \sigma^2 \quad \text{and} \quad \lim_{n \to +\infty} \frac{(2n)^{3-2\gamma}}{L(1/(2n))^2} \operatorname{Var} V_{2n}(X) = \sigma^2. \quad (48)$$

We must compute the term $\operatorname{Cov}(V_n(X), V_{2n}(X))$. We have

$$\Delta X_k^{(n)} = \Delta X_{2k+1}^{(2n)} + \Delta X_{2k-1}^{(2n)} + 2\Delta X_{2k}^{(2n)}.$$

Therefore,

$$V_n(X) = \sum_{k=1}^{n-1} [(\Delta X_{2k+1}^{(2n)})^2 + (\Delta X_{2k-1}^{(2n)})^2 + 4(\Delta X_{2k}^{(2n)})^2$$
$$+ 2\Delta X_{2k+1}^{(2n)} \Delta X_{2k-1}^{(2n)} + 4\Delta X_{2k+1}^{(2n)} \Delta X_{2k}^{(2n)} + 4\Delta X_{2k-1}^{(2n)} \Delta X_{2k}^{(2n)}].$$

Moreover,

$$V_{2n}(X) = \sum_{j=1}^{2n-1} (\Delta X_j^{(2n)})^2.$$

To simplify, we will set the notation $\Delta X_l$ for $\Delta X_l^{(2n)}$ and $d_{lp}$ for $d_{lp}^{(2n)}$. To compute $\mathbb{E}[V_n(X)V_{2n}(X)]$, we use the Isserlis formulas (see Isserlis [14]), which yield

$$\mathrm{Cov}[(\Delta X_l)^2, (\Delta X_j)^2] = 2(\mathbb{E}[\Delta X_l \Delta X_j])^2 = 2d_{lj}^2$$

and

$$\mathrm{Cov}[\Delta X_l \Delta X_p, (\Delta X_j)^2] = 2\mathbb{E}[\Delta X_l \Delta X_j]\mathbb{E}[\Delta X_p \Delta X_j] = 2d_{lj}d_{pj}.$$

So we can check that

$$\mathrm{Cov}(V_n(X), V_{2n}(X)) = \sum_{i=1}^{6} S_i,$$

with

$$S_1 = 2\sum_{k=1}^{n-1} \sum_{j=1}^{2n-1} d_{2k+1,j}^2, \qquad S_2 = 2\sum_{k=1}^{n-1} \sum_{j=1}^{2n-1} d_{2k-1,j}^2,$$

$$S_3 = 8\sum_{k=1}^{n-1} \sum_{j=1}^{2n-1} d_{2k,j}^2, \qquad S_4 = 4\sum_{k=1}^{n-1} \sum_{j=1}^{2n-1} d_{2k-1,j}d_{2k+1,j},$$

$$S_5 = 8\sum_{k=1}^{n-1} \sum_{j=1}^{2n-1} d_{2k+1,j}d_{2k,j}, \qquad S_6 = 8\sum_{k=1}^{n-1} \sum_{j=1}^{2n-1} d_{2k-1,j}d_{2k,j}.$$

However, using the same techniques as in the Proof of Theorem 2, we obtain the formulas

$$\lim_{n \to +\infty} \frac{2^{1-\gamma}n^{3-2\gamma}}{L(1/n)L(1/(2n))} S_1 = 2^{\gamma-3}\sigma^2,$$

$$\lim_{n \to +\infty} \frac{2^{1-\gamma}n^{3-2\gamma}}{L(1/n)L(1/(2n))} S_4 = 2^{\gamma-2}\sigma_{1,\mathrm{cov}}^2,$$

$$\lim_{n \to +\infty} \frac{2^{1-\gamma}n^{3-2\gamma}}{L(1/n)L(1/(2n))} S_5 = 2^{\gamma-1}\sigma_{2,\mathrm{cov}}^2$$

and, likewise,

$$\lim_{n \to +\infty} \frac{2^{1-\gamma} n^{3-2\gamma}}{L(1/n)L(1/(2n))} S_2 = 2^{\gamma-3} \sigma^2,$$

$$\lim_{n \to +\infty} \frac{2^{1-\gamma} n^{3-2\gamma}}{L(1/n)L(1/(2n))} S_3 = 2^{\gamma-1} \sigma^2,$$

$$\lim_{n \to +\infty} \frac{2^{1-\gamma} n^{3-2\gamma}}{L(1/n)L(1/(2n))} S_6 = 2^{\gamma-1} \sigma^2_{2,\mathrm{cov}}.$$

Hence,

$$\lim_{n \to +\infty} \frac{2^{1-\gamma} n^{3-2\gamma}}{L(1/n)L(1/(2n))} \mathrm{Cov}(V_n(X), V_{2n}(X)) = 2^{\gamma-2} \sigma^2_*. \tag{49}$$

Therefore, (47) is a consequence of (48) and (49).

Now we apply the Lindeberg central limit theorem to $S_n(\lambda, \mu)$ in the same manner as in Theorem 2. We set

$$\widetilde{S_n}(\lambda, \mu) = S_n(\lambda, \mu) - \mathbb{E} S_n(\lambda, \mu).$$

Because $\lambda \geq 0$ and $\mu \geq 0$, we can consider $S_n(\lambda, \mu)$ as the Euclidean norm of the Gaussian vector $(G_i; 1 \leq i \leq 3n - 2)$:

$$G_i = \sqrt{\lambda \frac{n^{1-\gamma}}{L(1/n)}} \Delta X_i^{(n)}, \qquad\qquad 1 \leq i \leq n-1,$$

$$G_i = \sqrt{\mu \frac{(2n)^{1-\gamma}}{L(1/(2n))}} \Delta X_{i+1-n}^{(2n)}, \qquad n \leq i \leq 3n-2.$$

Therefore, Cochran's theorem yields

$$S_n(\lambda, \mu) = \sum_{j=1}^{a_n} \tau_{j,n} (Y_n^{(j)})^2,$$

with $a_n$, $\tau_{j,n}$, $\tau_n^*$ and $Y_n^{(j)}$ as in the Proof of Theorem 2. This yields

$$\widetilde{S_n}(\lambda, \mu) = \sum_{j=1}^{a_n} \tau_{j,n} [(Y_n^{(j)})^2 - 1].$$

Also, we notice that

$$\tau_n^* \leq \max_{1 \leq j \leq 3n-2} \sum_{i=1}^{3n-2} |\mathbb{E}(Z_i Z_j)|.$$

Moreover, if $1 \le i \le n-1$ and $1 \le j \le 2n-1$,

$$\mathbb{E}(\Delta X_i^{(n)} \Delta X_j^{(2n)}) = \mathbb{E}[(\Delta X_{2i+1}^{(2n)} + \Delta X_{2i-1}^{(2n)} + 2\Delta X_{2i}^{(2n)})\Delta X_j^{(2n)}]$$
$$= d_{2i+1,j}^{(2n)} + d_{2i-1,j}^{(2n)} + 2d_{2i,j}^{(2n)},$$

so

$$\begin{aligned}
\tau_n^* \le K\, \frac{n^{1-\gamma}}{L(1/n)} \Bigg( & \lambda \max_{1 \le j \le n-1} \sum_{i=1}^{n-1} |d_{ij}^{(n)}| + \mu \max_{1 \le j \le 2n-1} \sum_{i=1}^{2n-1} |d_{ij}^{(2n)}| \\
& + \sqrt{\lambda\mu} \max_{1 \le j \le n-1} \sum_{i=1}^{2n-1} (|d_{i,2j+1}^{(2n)}| + |d_{i,2j-1}^{(2n)}| + 2|d_{i,2j}^{(2n)}|) \\
& + \sqrt{\lambda\mu} \max_{1 \le j \le 2n-1} \sum_{i=1}^{n-1} (|d_{2i+1,j}^{(2n)}| + |d_{2i-1,j}^{(2n)}| + 2|d_{2i,j}^{(2n)}|) \Bigg) \\
\le K\, \frac{n^{1-\gamma}}{L(1/n)} \Bigg( & \max_{1 \le j \le n-1} \sum_{i=1}^{n-1} |d_{ij}^{(n)}| + \max_{1 \le j \le 2n-1} \sum_{i=1}^{2n-1} |d_{ij}^{(2n)}| \Bigg).
\end{aligned}$$

Therefore, with the arguments of Bégyn [5], we can checked that

$$\tau_n^* \le K n^{\gamma-2} L\left(\frac{1}{n}\right).$$

Thus

$$0 \le \frac{\tau_n^*}{\sqrt{\mathrm{Var}\, S_n(\lambda,\mu)}} \le \frac{K}{\sqrt{n}}.$$

Using the Lindeberg central limit theorem, we obtain

$$\frac{S_n(\lambda,\mu) - \mathbb{E}S_n(\lambda,\mu)}{\sqrt{\mathrm{Var}\, S_n(\lambda,\mu)}} \xrightarrow{(\mathcal{L})} \mathcal{N}(0,1).$$

Hence, as in the Proof of Theorem 2, we establish the convergence announced in (46). With a generalization of the Cramér–Wold arguments, based on the properties of the Laplace transform that can be found in Istas and Lang ([15], page 431), we ascertain that the Laplace transform of the vector

$$\sqrt{n} \begin{pmatrix} \frac{n^{1-\gamma}}{L(\frac{1}{n})} V_n(X) - \int_0^1 g_0(t)\,\mathrm{d}t - \sum_{i=1}^q \int_0^1 g_i(t)\,\mathrm{d}t \cdot \phi(\frac{1}{n})^{\nu_i} \\ \frac{(2n)^{1-\gamma}}{L(\frac{1}{2n})} V_{2n}(X) - \int_0^1 g_0(t)\,\mathrm{d}t - \sum_{i=1}^q \int_0^1 g_i(t)\,\mathrm{d}t \cdot \phi(\frac{1}{2n})^{\nu_i} \end{pmatrix}$$

tends on $\mathbb{R}^2$ to the Laplace transform of a Gaussian centered law with covariance matrix $\Sigma$. This proves the result. $\qquad\square$

As we stated in the Introduction, we want to apply these results to estimate the parameters of some fractional processes. To explain how to use these results, we apply them to the FBM, even though the consequences are not new in this case. We refer to Sections 3 and 4 for original consequences.

## 2.4. Application to fractional Brownian motion

We study the example of the FBM $B^H$. We can check that the theorems of Section 2 can be applied with $\gamma = 2 - 2H$, $L(h) = 1$, $q = 0$, $g_0(t) = 4 - 2^{2H}$ and $\phi(h) = h$.

By applying Theorem 1, it follows that almost surely (see Cohen *et al.* [10])

$$\lim_{n \to +\infty} n^{2H-1} V_n(B^H) = 4 - 2^{2H}. \tag{50}$$

Next, if we apply Theorem 2, then we obtain (see Coeurjolly [9])

$$\sqrt{n}(n^{2H-1} V_n(B^H) - (4 - 2^{2H})) \xrightarrow{(\mathcal{L})} \mathcal{N}(0, \sigma^2_{\text{FBM},H}) \tag{51}$$

with

$$\sigma^2_{\text{FBM},H} = 2(4 - 2^{2H})^2 + (2^{2H+2} - 7 - 3^{2H})^2$$
$$+ (2H)^2(2H-1)^2(2H-2)^2(2H-3)^2 \|\rho_{2-2H}\|^2, \tag{52}$$

because, when computed, $C(s,t) = -H(2H-1)(2H-2)(2H-3)$ and $\widetilde{g}(t) = (2^{2H+2} - 7 - 3^{2H})/2$.

Because of Theorem 3, we get

$$\sqrt{n} \begin{pmatrix} n^{2H-1} V_n(B^H) - (4 - 2^{2H}) \\ (2n)^{2H-1} V_{2n}(X) - (4 - 2^{2H}) \end{pmatrix} \xrightarrow{(\mathcal{L})} \mathcal{N}(0, \Sigma_{\text{FBM},H}) \tag{53}$$

with

$$\Sigma_{\text{FBM},H} = \begin{pmatrix} \sigma^2_{\text{FBM},H} & 2^{-2H} \sigma^2_{*,\text{FBM},H} \\ 2^{-2H} \sigma^2_{*,\text{FBM},H} & 2^{-1} \sigma^2_{\text{FBM},H} \end{pmatrix} \tag{54}$$

and

$$\sigma^2_{1,\text{cov},\text{FBM},H} = (2H)^2(2H-1)^2(2H-2)^2(2H-3)^2 \sum_{l=2}^{+\infty} \rho_{2-2H}(l) \rho_{2-2H}(l-2)$$
$$+ \frac{1}{2}(2^{2H+2} - 7 - 3^{2H})^2, \tag{55}$$

$$\sigma^2_{2,\text{cov},\text{FBM},H} = (2H)^2(2H-1)^2(2H-2)^2(2H-3)^2 \sum_{l=2}^{+\infty} \rho_{2-2H}(l) \rho_{2-2H}(l-1)$$
$$+ 2(4 - 2^{2H})(2^{2H+2} - 7 - 3^{2H}), \tag{56}$$

$$\sigma^2_{*,\text{FBM},H} = 3\sigma^2_{\text{FBM},H} + \sigma^2_{1,\text{cov},\text{FBM},H} + 4\sigma^2_{2,\text{cov},\text{FBM},H}. \tag{57}$$

The $\delta$ method yields that the statistic

$$\widehat{H_n} = \frac{1}{2} - \log \frac{V_{2n}(B^H)}{V_n(B^H)} (2 \log 2)^{-1}$$

is a strongly consistent estimator of $H$ and that (see Coeurjolly [9])

$$\sqrt{n}(\widehat{H_n} - H) \overset{(\mathcal{L})}{\longrightarrow} \mathcal{N}\left(0, \frac{3\sigma_{\text{FBM},H}^2 - 2^{2-2H}\sigma_{*,\text{FBM},H}^2}{4(4 - 2^{2H})\log 2}\right).$$

# 3. Two-parameter fractional Brownian motion

Two-parameter fractional Brownian motion was introduced by Houdré and Villa [13] as an example of a quasi-helix. *Two-parameter fractional Brownian motion* $B^{H,K} = \{B_t^{H,K}; t \in \mathbb{R}\}$ is defined, for $H \in ]0, 1[$ and $K \in ]0, 1]$, as the unique continuous centered Gaussian process with covariance function

$$\forall s, t \in \mathbb{R} \qquad R^{H,K}(s, t) = \text{Cov}(B_s^{H,K}, B_t^{H,K}) = \frac{1}{2^K}((s^{2H} + t^{2H})^K - |s - t|^{2HK}).$$

The process $B^{H,K}$ is $HK$-self-similar, it has a critical Hölder exponent equal to $HK$ in the sense of Adler [1], it is, for $K = 1$, the standard fractional Brownian motion and it has stationary increments if and only if $K = 1$.

We refer to Houdré and Villa [13] for the proofs. In their paper, they introduced the process $B^{H,K}$ under the name *bifractional Brownian motion*. We suggest here to call it *two-parameter fractional Brownian motion*, because it is a monofractal process (the pointwise Hölder exponent of its trajectories is a.s. constant) and the term 'bifractional' may imply that it is a multifractal process with two values for its pointwise Hölder exponent.

Theorems 1 and 2 yield estimators of $H$ and $K$, the relevant quantities in the study of $B^{H,K}$. To determine the estimators, that we assume that we dispose of the observation of one path of $B^{H,K}$ on the interval $[T_1, T_2]$, where $T_1, T_2 \in \mathbb{R}$, $T_1 < T_2$. This process is considered to be indexed by $[0, 1]$. Therefore, we introduce the process $Y$ defined by

$$\forall t \in [0, 1] \qquad Y_t = B^{H,K}(\tau(t)) \qquad \text{with } \tau(t) = (T_2 - T_1)t + T_1.$$

We obtain a new process $Y$, which is centered, Gaussian and has covariance function $r^{H,K}(s, t) = \text{Cov}(Y_s, Y_t) = R^{H,K}(\tau(s), \tau(t))$, and we dispose of the observation of one path of $Y$ on $[0, 1]$. The results of Section 2 can be applied to $Y$ under the condition $[T_1, T_2] \subset ]0, +\infty[$.

## 3.1. The results

First we study the almost sure convergence of the second-order quadratic variations.

**Proposition 4.** *We have when* $n \to +\infty$,

$$n^{2HK-1} V_n(Y) \xrightarrow{a.s.} \frac{4 - 2^{2HK}}{2^{K-1}} (T_2 - T_1)^{2HK}. \tag{58}$$

Next we study the weak convergence.

**Proposition 5.** *We have, when* $n \to +\infty$,

$$\sqrt{n} \left( n^{2HK-1} V_n(Y) - \frac{4 - 2^{2HK}}{2^{K-1}} (T_2 - T_1)^{2HK} \right)$$

$$\xrightarrow{(\mathcal{L})} \mathcal{N} \left( 0, \frac{(T_2 - T_1)^{4HK}}{2^{2(K-1)}} \sigma^2_{\mathrm{FBM}, HK} \right), \tag{59}$$

*where* $\sigma^2_{\mathrm{FBM}, HK}$ *was defined in* (52).

As in the case of the FBM, we can deduce an estimator of $HK$.

**Proposition 6.** *The statistic*

$$\widehat{HK}_n = \frac{1}{2} - \log \frac{V_{2n}(Y)}{V_n(Y)} (2 \log 2)^{-1} \tag{60}$$

*is a strongly consistent estimator of* $HK$ *and when* $n \to +\infty$,

$$\sqrt{n}(\widehat{HK}_n - HK) \xrightarrow{(\mathcal{L})} \mathcal{N} \left( 0, \frac{3\sigma^2_{\mathrm{FBM}, HK} - 2^{2-2HK} \sigma^2_{*, \mathrm{FBM}, HK}}{2^{K+1}(4 - 2^{2HK}) \log 2} (T_2 - T_1)^{4HK} \right), \tag{61}$$

*where* $\sigma^2_{\mathrm{FBM}, HK}$ *and* $\sigma^2_{*, \mathrm{FBM}, HK}$ *were defined in* (52) *and* (57).

The quantity $HK$ is relevant in the study of $B^{H,K}$, but it does not characterize the law of this process. For this characterization, we need to know both parameters $H$ and $K$. A refinement of the previous results enables us to construct strongly consistent and asymptotically normal estimators of these quantities.

**Proposition 7.** *The statistic*

$$\widehat{K}_n = 1 - \frac{1}{\log 2} \log \left( \frac{n^{2\widehat{HK}_n - 1}}{(4 - 2^{2\widehat{HK}_n})(T_2 - T_1)^{2\widehat{HK}_n}} V_n(Y) \right) \tag{62}$$

*is a strongly consistent estimator of* $K$ *and when* $n \to +\infty$,

$$\sqrt{n}(\widehat{K}_n - K) \xrightarrow{(\mathcal{L})} \mathcal{N} \left( 0, \frac{\sigma^2_{\mathrm{FBM}, HK}}{(4 - 2^{2HK})^2 \log^2 2} \right), \tag{63}$$

*where $\sigma^2_{\mathrm{FBM},HK}$ was defined in (52). Moreover the statistic*

$$\widehat{H}_n = \frac{\widehat{HK}_n}{\widehat{K}_n} \tag{64}$$

*is a strongly consistent estimator of $H$ and when $n \to +\infty$,*

$$\sqrt{n}(\widehat{H}_n - H) \xrightarrow{(\mathcal{L})} \mathcal{N}\left(0, \frac{(T_2 - T_1)^{4HK}}{2^{2(K-1)}}\eta^2\right) \tag{65}$$

*with*

$$\eta^2 = \frac{H^2}{K^2}\eta_1 + \frac{1}{K^2}\eta_2 - \frac{2H}{K^2}\eta_3$$

*and*

$$\eta_1 = \frac{2^{2K-2}\sigma^2_{\mathrm{FBM},HK}}{(4 - 2^{2HK})^2 \log^2 2 (T_2 - T_1)^{4HK}},$$

$$\eta_2 = 2^{K-1}\frac{3\sigma^2_{\mathrm{FBM},HK} - 2^{2-2HK}\sigma^2_{*,\mathrm{FBM},HK}}{4(4 - 2^{2HK})\log 2},$$

$$\eta_3 = 2^{2K-2}\frac{2^{-2HK}\sigma^2_{*,\mathrm{FBM},HK} - \sigma^2_{\mathrm{FBM},HK}}{2(4 - 2^{2HK})^2 \log^2 2 (T_2 - T_1)^{2HK}},$$

*where $\sigma^2_{*,\mathrm{FBM},HK}$ was defined in (57).*

## 3.2. Proofs of the results for two-parameter FBM

**Proof of Proposition 4.** We apply Theorem 1 to the process $Y$. We need to show only that assumptions 2(b) and 2(c) (in Theorem 1) are satisfied (the other assumption is obvious).

For assumption 2(b), it is clear that the derivative $\frac{\partial^4 r^{H,K}}{\partial s^2 \partial t^2}(s, t)$ exists on $]0, 1]^2 \setminus \{s = t\}$. Moreover, we can check that, $\forall s, t \in ]0, 1]^2 \setminus \{s = t\}$,

$$\frac{\partial^4 r^{H,K}}{\partial s^2 \partial t^2}(s, t) = -\frac{2HK(2HK-1)(2HK-2)(2HK-3)}{2^K}(T_2 - T_1)^{2HK}|s - t|^{2HK-4}$$
$$+ (T_2 - T_1)^4 \psi(\tau(s), \tau(t)), \tag{66}$$

where $\psi(\tau(s), \tau(t))$ is continuous on $[0, 1]^2$. Therefore, the assumption 2(b) (in Theorem 1) is satisfied with $L(h) = 1$ and $\gamma = 2 - 2HK$.

For assumption 2(c) (in Theorem 1), computations yield

$$\frac{(\delta_1^h \circ \delta_2^h r^{H,K})(t, t)}{h^{2HK}} = \frac{4 - 2^{2HK}}{2^{K-1}}(T_2 - T_1)^{2HK} + \frac{\varepsilon_t(h)}{h^{2HK}} \tag{67}$$

and we can check that $\varepsilon_t(0) = \varepsilon'_t(0) = \varepsilon''_t(0) = \varepsilon_t^{(3)}(0) = 0$. So that Taylor formula yields

$$\forall h \le t \le 1 - h \qquad \varepsilon_t(h) = \int_0^h \frac{(h-x)^3}{3!} \varepsilon_t^{(4)}(x) \, \mathrm{d}x.$$

Therefore, we have

$$\sup_{h \le t \le 1-h} \sup_{0 \le x \le h} |\varepsilon_t^{(4)}(x)| \stackrel{h \to 0^+}{=} \mathcal{O}(1),$$

which yields

$$\sup_{h \le t \le 1-h} \left| \frac{(\delta_1^h \circ \delta_2^h r^{H,K})(t,t)}{h^{2HK}} - (T_2 - T_1)^{2HK} \frac{4 - 2^{2HK}}{2^{K-1}} \right| \stackrel{h \to 0^+}{=} \mathcal{O}(h^{4-2HK}). \qquad (68)$$

Therefore, the assumption 2(c) (in Theorem 1) is fulfilled with

$$g_0(t) = \frac{4 - 2^{2HK}}{2^{K-1}} (T_2 - T_1)^{2HK}. \qquad (69)$$

Consequently, we can appply Theorem 1 to $Y$ and obtain (58). □

**Proof of Proposition 5.** We apply Theorem 2 to the process $Y$. As in the Proof of Proposition 4, we need to show only that the assumptions 2 and 3 (in Theorem 2) are satisfied.

For assumption 2 (in Theorem 2) the previous proof showed formula (66), which yields that for all $s, t \in ]0,1]^2 \setminus \{s = t\}$,

$$\frac{\partial^4 r^{H,K}}{\partial s^2 \, \partial t^2}(s,t) = -\frac{2HK(2HK-1)(2HK-2)(2HK-3)}{2^K} (T_2 - T_1)^{2HK} |s-t|^{2HK-4}$$
$$+ (T_2 - T_1)^4 \psi(\tau(s), \tau(t)),$$

where $\psi(\tau(s), \tau(t))$ is continuous on $]0,1]^2$. Therefore, the assumption (21) in Theorem 2 is satisfied with $L(h) = 1$, $\gamma = 2 - 2HK$ and

$$C(s,t) = -\frac{2HK(2HK-1)(2HK-2)(2HK-3)}{2^K} (T_2 - T_1)^{2HK}$$
$$+ (T_2 - T_1)^4 |s-t|^{4-2HK} \psi(\tau(s), \tau(t)). \qquad (70)$$

For assumption 3 (in Theorem 2), formula (68) of the previous proof shows that the assumption 3(d) (in Theorem 2) is fulfilled with $q = 0$, $g_0(t) = \frac{8 - 2^{2HK+1}}{2^K}$ and $\alpha_0 = 1/2$. Moreover, we can check that

$$\frac{(\delta_1^h \circ \delta_2^h r^{H,K})(t, t+h)}{h^{2HK}} = \frac{2^{2HK+2} - 3^{2HK} - 7}{2^K} (T_2 - T_1)^{2HK} + \frac{\eta_t(h)}{h^{2HK}}.$$

With the same arguments as those used for $\varepsilon_t(h)$ in the previous proof, we obtain

$$\sup_{h \le t \le 1-h} |\eta_t(h)| \stackrel{h \to 0^+}{=} \mathcal{O}(h^4).$$

This shows that the assumption 3(e) (in Theorem 2) is satisfied with

$$\widetilde{g}(t) = \frac{2^{2HK+2} - 3^{2HK} - 7}{2^K}(T_2 - T_1)^{2HK}. \tag{71}$$

Consequently, we can appply Theorem 2 to $Y$ and obtain (59). $\qquad \square$

**Proof of Proposition 6.** We apply the $\delta$ method with the $C^1$ function

$$f(x, y) = \frac{1}{2} - \frac{\log(y/x)}{2\log 2}$$

to the convergence announced in (44) to yield the result. $\qquad \square$

**Proof of Proposition 7.** First we establish a refinement of Proposition 4. Because of (68), we have, for all $\alpha \in {]}0, 1{[}$,

$$\sup_{h \le t \le 1-h} \left| \frac{(\delta_1^h \circ \delta_2^h r^{H,K})(s,t)}{h^{2HK}} - (T_2 - T_1)^{2HK} \frac{8 - 2^{2HK+1}}{2^K} \right| \stackrel{h \to 0^+}{=} o(h^\alpha).$$

Therefore, the assumption 2(c) of Theorem 1 is fulfilled with $q = 1$, $g_0(t) = (T_2 - T_1)^{2HK} \times \frac{4 - 2^{2HK}}{2^{K-1}}$, $g_1(t) = 0$, $\phi(h) = h$ and $\nu_1 = \alpha$. It yields that almost surely

$$n^{2HK-1} V_n(Y) \stackrel{n \to +\infty}{=} (T_2 - T_1)^{2HK} \frac{4 - 2^{2HK}}{2^{K-1}} + o\left(\frac{1}{n^\alpha}\right).$$

Taylor expansions yield that almost surely

$$\widehat{HK}_n \stackrel{n \to +\infty}{=} HK + o\left(\frac{1}{n^\alpha}\right)$$

and

$$\frac{n^{2\widehat{HK}_n - 1}}{n^{2HK-1}} \stackrel{n \to +\infty}{=} 1 + o\left(\frac{\log n}{n^\alpha}\right).$$

With $\alpha = 3/4$, we obtain that almost surely

$$\frac{n^{2\widehat{HK}_n - 1}}{n^{2HK-1}(4 - 2^{2\widehat{HK}_n})(T_2 - T_1)^{2\widehat{HK}_n}} \stackrel{n \to +\infty}{=} \frac{1}{(4 - 2^{2HK})(T_2 - T_1)^{2HK}} + o\left(\frac{\log n}{n^{3/4}}\right). \tag{72}$$

In addition, we have (44):

$$\sqrt{n}\begin{pmatrix} n^{2HK-1}V_n(Y) - \frac{4-2^{2HK}}{2^{K-1}}(T_2-T_1)^{2HK} \\ (2n)^{2HK-1}V_{2n}(Y) - \frac{4-2^{2HK}}{2^{K-1}}(T_2-T_1)^{2HK} \end{pmatrix} \xrightarrow{(\mathcal{L})} \mathcal{N}\left(0, \frac{(T_2-T_1)^{4HK}}{2^{2(K-1)}}\Sigma_{\mathrm{FBM},HK}\right).$$

If we apply the $\delta$ method with the $C^1$ function

$$f(x,y) = \begin{pmatrix} x \\ \frac{1}{2} - \frac{\log(y/x)}{2\log 2} \end{pmatrix},$$

then the Slutsky lemma and (72) yield that there exists a $2 \times 2$ real matrix $\mathbb{A}$ such that

$$\sqrt{n}\begin{pmatrix} \frac{n^{2\widehat{HK}_n-1}}{(4-2^{2\widehat{HK}_n})(T_2-T_1)^{2\widehat{HK}_n}}V_n(Y) - \frac{1}{2^{K-1}} \\ \widehat{HK}_n - HK \end{pmatrix} \xrightarrow{(\mathcal{L})} \mathcal{N}\left(0, \frac{(T_2-T_2)^{4HK}}{2^{2(K-1)}}\mathbb{A}\right).$$

By again applying the $\delta$ method with the $C^1$ function

$$f(x,y) = \begin{pmatrix} 1 - \frac{\log x}{\log 2} \\ y \end{pmatrix},$$

we obtain

$$\sqrt{n}\begin{pmatrix} \widehat{K}_n - K \\ \widehat{HK}_n - HK \end{pmatrix} \xrightarrow{(\mathcal{L})} \mathcal{N}\left(0, \frac{(T_2-T_2)^{4HK}}{2^{2(K-1)}}\mathbb{C}\right)$$

with

$$\mathbb{C} = \begin{pmatrix} \eta_1 & \eta_3 \\ \eta_3 & \eta_2 \end{pmatrix}.$$

This proves (63). A final application of the $\delta$ method with the $C^1$ function

$$f(x,y) = \frac{y}{x}$$

yields (65). □

# 4. Anisotropic fractional Brownian motion

Let $d \in \mathbb{N}^*$. Let $H : \mathbb{R}^d \to ]0,1[$ be a Borelian function that is homogeneous of degree zero,

$$\forall \xi \in \mathbb{R}^d, \forall \lambda \in \mathbb{R} \setminus \{0\} \qquad H(\lambda\xi) = H(\xi),$$

that can be identified with an even function from the sphere $S^{d-1}$ into $\mathbb{R}$ that we denote $H$ as well. We assume, moreover, that $H$ takes its values inside the interval $[\underline{H}, \overline{H}] \subset ]0,1[$, with $\underline{H} = \mathrm{ess\,inf}\,H$ and $\overline{H} = \mathrm{ess\,sup}\,H$.

We define the anisotropic fractional Brownian motion (AFBM) with directional Hurst index $H$, denoted $A^{(H)}$, by the harmonizable representation formula

$$\forall u \in \mathbb{R}^d \qquad A^{(H)}(u) = \int_{\mathbb{R}^d} \frac{e^{i\langle u, \xi \rangle} - 1}{|\xi|^{H(\xi) + d/2}} \, dW(\xi), \tag{73}$$

where $\langle \cdot, \cdot \rangle$ is the canonical scalar product and $W$ is a complex random measure in the sense of Samorodnitsky and Taqqu ([16], 325–328). It is a Gaussian field with stationary increments. Bonami and Estrade [8] showed that $A^{(H)}$ has a critical Hölder exponent equal to $\underline{H}$. Moreover, they showed that the field $A^{(H)}$ is locally asymtotically self-similar (l.a.s.s.) of order $\underline{H}$ at any point of $\mathbb{R}^d$ (see Definition 8) if and only if $\text{Leb}(\{H(\theta) = \underline{H}\}) > 0$, where $\text{Leb}(\cdot)$ denotes indifferently the Lebesgue measure on $\mathbb{R}^d$ or the Lebesgue measure on $S^{d-1}$.

Let us recall the definition of the l.a.s.s. property (see Benassi, Jaffard and Roux [4]):

**Definition 8.** *Let $\beta > 0$. A process $\{X_u; u \in \mathbb{R}^d\}$ is locally asymptotically self-similar (l.a.s.s.) of order $\beta$ at point $u_0 \in \mathbb{R}^d$ if the finite-dimensional distributions of the process*

$$\left\{ \frac{X(u_0 + \lambda u) - X(u_0)}{\lambda^{\beta}}; u \in \mathbb{R}^d \right\}$$

*converge to the finite-dimensional distributions of a non-zero Gaussian process when $\lambda \to 0^+$. The limit process is called the tangent process at point $u_0$.*

Our purpose is to identify the function $H$ when we consider one realization of the field $A^{(H)}$. For that we apply the theorems shown in Section 2 and restrict the field to some segment of $\mathbb{R}^d$.

To simplify the computations we assume next that $d = 2$. Note that in this case we can identify $H$ with an even $\pi$-periodic function on $\mathbb{R}$. We consider one realization of $A^{(H)}$, which is observed in axes denoted by $Oxy$. We assume too that these axes of observation are equal to the axes of definition of $A^{(H)}$.

Let $[A, B]$ be the radial segment of length $L \in ]0, +\infty[$ such that the distance between $O$ and $A$ is equal to $L\varepsilon$ (with $\varepsilon \geq 0$) and the angle between $[A, B]$ and the axes $Ox$ is equal to $\omega \in [0, 2\pi[$. See Figure 1 for more details on the geometry of the problem (and note that the angle are oriented anticlockwise).

We use the following parametrization of the point $u = (u_1, u_2) \in [A, B]$:

$$\begin{aligned} u_1 &= L(t + \varepsilon) \cos \omega, \\ u_2 &= L(t + \varepsilon) \sin \omega, \end{aligned} \tag{74}$$

where $t$ goes over the interval $[0, 1]$. Next we consider the restriction of the field $A^{(H)}$ to the segment $[A, B]$ with the parametrization (74). Hence we obtain a new process $Z$ indexed by $t \in [0, 1]$:

$$Z_t = A^{(H)} \begin{pmatrix} L(t + \varepsilon) \cos \omega \\ L(t + \varepsilon) \sin \omega \end{pmatrix}.$$

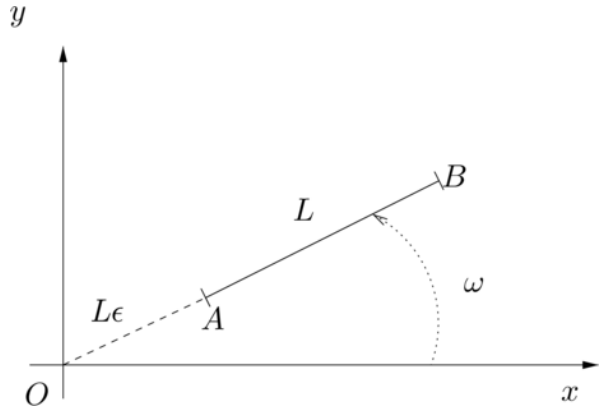

**Figure 1.** Geometry of the problem.

It is clear that $Z$ is a centered Gaussian process with stationary increments.

To apply the theorems of Section 2, we must compute the covariance function $R$ of $Z$ that is given by the following lemma.

**Lemma 9.** *The covariance function of the process* $(Z_t)_{t \in [0,1]}$ *is given by, for all* $s, t \in [0,1]$,

$$R(s,t) = 4 \int_0^\pi \Lambda(\theta)[|s + \varepsilon|^{2H(\theta)} + |t + \varepsilon|^{2H(\theta)} - |s - t|^{2H(\theta)}] \,\mathrm{d}\theta, \tag{75}$$

*with*

$$\forall \theta \in [0, 2\pi] \qquad \Lambda(\theta) = \frac{C(1, H(\theta))}{8C(2, \underline{H})^2} L^{2H(\theta)} |\cos(\theta - \omega)|^{2H(\theta)}, \tag{76}$$

*where for all* $d \in \mathbb{N}^*$ *and for all* $H \in ]0, 1[$,

$$C(d, H) = \left( \frac{\pi^{(d+1)/2} \Gamma(H + 1/2)}{H \Gamma(2H) \sin(H\pi) \Gamma(H + d/2)} \right)^{1/2} \tag{77}$$

*and* $\Gamma$ *denotes the Euler gamma function.*

**Proof of Lemma 9.** First we compute the variogram $v$ of the field $A^{(H)}$. For all $u \in \mathbb{R}^2$, we have

$$v(u) = \frac{1}{4} \mathbb{E}[(A^{(H)}(u))^2] = \frac{1}{8} \int_{S^1} C(1, H(y))^2 |\langle u, y \rangle|^{2H(y)} \,\mathrm{d}y$$

$$= \frac{1}{8C(2, \underline{H})^2} \int_0^{2\pi} C(1, H(\theta)) |u_1 \cos(\theta) + u_2 \sin(\theta)|^{2H(\theta)} \,\mathrm{d}\theta.$$

Because the field $A^{(H)}$ has stationary increments and is vanishing a.s. at the origin, its covariance function is given by

$$\forall u, u' \in \mathbb{R}^2, \qquad \text{Cov}(A^{(H)}(u), A^{(H)}(u')) = 2(v(u) + v(u') - v(u - u')).$$

If we take $u, u'$ in the segment $[A, B]$ with $u$ parametrized by $s$ and $u'$ by $t$, we obtain

$$R(s, t) = \text{Cov}(A^{(H)}(u), A^{(H)}(u'))$$

$$= 2 \int_0^{2\pi} \Lambda(\theta)[|s + \varepsilon|^{2H(\theta)} + |t + \varepsilon|^{2H(\theta)} - |s - t|^{2H(\theta)}] \, d\theta.$$

Because the functions $\Lambda$ and $H$ are $\pi$-periodic, we obtain (75). $\qquad \square$

By applying theorems of Section 2 to the process $Z$, we are able to estimate $\underline{H}$, the Hölder critical exponent of $A^{(H)}$. We distinguish two cases.

## 4.1. The l.a.s.s. case

In this subsection we assume that

$$\text{Leb}\{H(\theta) = \text{ess inf } H\} > 0$$

and we use the notation

$$\forall \underline{H} \in \,]0, 1[ \qquad J_{\underline{H}} = 8 \int_0^\pi \Lambda(\theta) \mathbf{1}_{\{H(\theta) = \underline{H}\}} \, d\theta. \tag{78}$$

**Proposition 10.** *We have, when* $n \to +\infty$,

$$n^{2\underline{H}-1} V_n(Z) \xrightarrow{a.s.} (4 - 2^{2\underline{H}}) J_{\underline{H}}. \tag{79}$$

For the central limit theorem we must study the asymptotic of

$$h \mapsto 8 \int_0^\pi \Lambda(\theta)(4 - 2^{2H(\theta)})|h|^{2(H(\theta) - \underline{H})} \mathbf{1}_{\{\underline{H} < H(\theta) \le \underline{H} + 1/4\}} \, d\theta$$

when $h \to 0^+$. For this purpose, we need to sharpen the assumption on the function $H$. We get two kinds of central limit theorem.

*Case I.* Assume that

$$\text{Leb}(\{\underline{H} < H(\theta) \le \underline{H} + 1/4\}) = 0.$$

**Proposition 11.** *We have, when $n \to +\infty$,*

$$\sqrt{n}(n^{2\underline{H}-1}V_n(Z) - (4 - 2^{2\underline{H}})J_{\underline{H}}) \xrightarrow{(\mathcal{L})} \mathcal{N}(0, J_{\underline{H}}^2 \sigma_{\mathrm{FBM},\underline{H}}^2). \tag{80}$$

From these results, we can deduce a strongly consistent estimator of $\underline{H}$ that is asymptotically normal.

**Corollary 12.** *The statistic*

$$\widehat{\underline{H}}_n = \frac{1}{2} - \log \frac{V_{2n}(Z)}{V_n(Z)}(2\log 2)^{-1} \tag{81}$$

*is a strongly consistent estimator of $\underline{H}$ and*

$$\sqrt{n}(\widehat{\underline{H}}_n - \underline{H}) \xrightarrow{(\mathcal{L})} \mathcal{N}\left(0, J_{\underline{H}} \frac{3\sigma_{\mathrm{FBM},\underline{H}}^2 - 2^{2-2\underline{H}}\sigma_{*,\mathrm{FBM},\underline{H}}^2}{4(4 - 2^{2\underline{H}})\log 2}\right), \tag{82}$$

*where $\sigma_{\mathrm{FBM},\underline{H}}^2$ and $\sigma_{*,\mathrm{FBM},\underline{H}}^2$ are as defined in (52) and (57).*

*Case II.* Assume that

$$\mathrm{Leb}(\{\underline{H} < H(\theta) \leq \underline{H} + 1/4\}) > 0.$$

**Proposition 13.** *We have almost surely*

$$n^{2\underline{H}-1}V_n(Z) \overset{n \to +\infty}{\equiv} (4 - 2^{2\underline{H}})J_{\underline{H}} + \phi\left(\frac{1}{n}\right) + o\left(\phi\left(\frac{1}{n}\right)\right) \tag{83}$$

*with*

$$\phi(h) = 8\int_0^\pi \Lambda(\theta)(4 - 2^{2H(\theta)})|h|^{2(H(\theta)-\underline{H})}\mathbf{1}_{\{\underline{H} < H(\theta) \leq \underline{H}+1/4\}}\,\mathrm{d}\theta.$$

*Moreover, when $n \to +\infty$,*

$$\sqrt{n}\left(n^{2\underline{H}-1}V_n(Z) - (4 - 2^{2\underline{H}})J_{\underline{H}} - \phi\left(\frac{1}{n}\right)\right) \xrightarrow{(\mathcal{L})} \mathcal{N}(0, J_{\underline{H}}^2 \sigma_{\mathrm{FBM},\underline{H}}^2). \tag{84}$$

As in Case I, we can deduce a strongly consistent estimator of $\underline{H}$ that is asymptotically normal.

**Corollary 14.** *The statistic*

$$\widehat{\underline{H}}_n = \frac{1}{2} - \log \frac{V_{2n}(Z)}{V_n(Z)}(2\log 2)^{-1}, \tag{85}$$

*is a strongly consistent estimator of $\underline{H}$ and*

$$\sqrt{n}\left(\widehat{\underline{H}}_n - \underline{H} + \log\left(\frac{(4 - 2^{2\underline{H}})J_{\underline{H}} + \phi(1/(2n))}{(4 - 2^{2\underline{H}})J_{\underline{H}} + \phi(1/n)}\right)(2\log 2)^{-1}\right)$$

$$\overset{(\mathcal{L})}{\longrightarrow} \mathcal{N}\left(0, J_{\underline{H}}\frac{3\sigma^2_{\mathrm{FBM},\underline{H}} - 2^{2 - 2\underline{H}}\sigma^2_{*,\mathrm{FBM},\underline{H}}}{4(4 - 2^{2\underline{H}})\log 2}\right). \tag{86}$$

**Remark.** In general, we have

$$\lim_{n \to +\infty} \log\left(\frac{(4 - 2^{2\underline{H}})J_{\underline{H}} + \phi(1/(2n))}{(4 - 2^{2\underline{H}})J_{\underline{H}} + \phi(1/n)}\right)(2\log 2)^{-1} = 0,$$

but the limit

$$\lim_{n \to +\infty} \sqrt{n} \log\left(\frac{(4 - 2^{2\underline{H}})J_{\underline{H}} + \phi(1/(2n))}{(4 - 2^{2\underline{H}})J_{\underline{H}} + \phi(1/n)}\right)(2\log 2)^{-1}$$

does not exist in general.

**Example.** We assume that

$$H(\theta) = \underline{H}\mathbf{1}_{]0,\theta_0[}(\theta) + \overline{H}\mathbf{1}_{]\theta_0,\pi[}(\theta)$$

with $\underline{H} < \overline{H}$ and $0 < \theta_0 < \pi$.

If $\overline{H} > \underline{H} + 1/4$, this is Case I.

If $\overline{H} = \underline{H} + 1/4$, this is Case II. We can check that $\phi(h) = (4 - 2^{2\underline{H}+1/2})J_{\underline{H}+1/4}h^{1/2}$, which implies

$$\lim_{n \to +\infty} \sqrt{n} \log\left(\frac{(4 - 2^{2\underline{H}})J_{\underline{H}} + \phi(1/(2n))}{(4 - 2^{2\underline{H}})J_{\underline{H}} + \phi(1/n)}\right)(2\log 2)^{-1} = \frac{(4 - 2^{2\underline{H}+1/2})J_{\underline{H}+1/4}}{(4 - 2^{2\underline{H}})J_{\underline{H}}}\frac{1 - \sqrt{2}}{2\sqrt{2}\log 2}.$$

Therefore, the Slutsky lemma yields that $\sqrt{N_n}(\widehat{\underline{H}}_n - \underline{H})$ converges in law to a Gaussian random variable with mean equal to the opposite of the right-hand side of the previous formula.

If $\overline{H} < \underline{H} + 1/4$, then this is also Case II and we have $\phi(h) = (4 - 2^{2\overline{H}})J_{\overline{H}}h^{2(\overline{H} - \underline{H})}$, which implies

$$\lim_{n \to +\infty} \sqrt{n} \log\left(\frac{(4 - 2^{2\underline{H}})J_{\underline{H}} + \phi(1/(2n))}{(4 - 2^{2\underline{H}})J_{\underline{H}} + \phi(1/n)}\right)(2\log 2)^{-1} = +\infty.$$

In this case, the bias term becomes infinite, so in practice it cannot be ignored.

## 4.2. The non-l.a.s.s. case

In this subsection we assume that

$$\mathrm{Leb}(\{H(\theta) = \mathrm{ess\,inf}\,H\}) = 0.$$

Next, we add the following assumptions:

1. $H$ is of class $C^2$ on $\mathbb{R}$ and there exists a unique point $\theta_* \in [0, \pi[$ such that $H(\theta_*) = \underline{H} = \inf H$.
2. $H$ is decreasing on $]\theta_* - \pi/2, \theta_*[$ and $H$ is increasing on $]\theta_*, \theta_* + \pi/2[$. Moreover, $H''(\theta_*) > 0$.
3. The angle $\omega$ is such that $\omega \neq \theta_* + \pi/2 \pmod{\pi}$.

These assumptions are classical when the Laplace method is applied. It is well known that they can be weakened, but we skip the technical details (see Dieudonné [11], page 125).

**Proposition 15.** *For all $t_0 \in \mathbb{R}^2$, the finite-dimensional distributions of*

$$\left\{ (-\log \varepsilon)^{1/4} \frac{A^{(H)}(t_0 + \varepsilon t) - A^{(H)}(t_0)}{\varepsilon^{\underline{H}}}; t \in \mathbb{R}^2 \right\} \tag{87}$$

*converge when $\varepsilon \to 0^+$ to the finite-dimensional distributions of the process*

$$\left\{ \frac{C(1, \underline{H})}{2} \left( \frac{1}{H''(\theta_*)} \right)^{1/4} B^{\underline{H}}(p_{\theta_*}(t)); t \in \mathbb{R}^2 \right\}, \tag{88}$$

*where $p_{\theta_*}(t)$ is the orthogonal projection of $t$ on the straight line $\{r\mathrm{e}^{\mathrm{i}\theta_*}, r \in \mathbb{R}\}$ and $B^{\underline{H}}$ is a standard FBM of Hurst index $\underline{H}$.*

Consequently, the field $A^{(H)}$ is no longer l.a.s.s. and we have shown that a normalization with a logarithm factor yields a non-trivial limit field. We will see in the sequel that we have the same behavior for the second-order quadratic variations.

For this section we use the notation

$$G_{\theta_*} = 8\Lambda(\theta_*)\sqrt{\frac{\pi}{H''(\theta_*)}}.$$

**Proposition 16.** *We have, when $n \to +\infty$,*

$$n^{2\underline{H}-1}\sqrt{\log n}\,V_n(Z) \xrightarrow{a.s.} (4 - 2^{2\underline{H}})G_{\theta_*}. \tag{89}$$

Now we want to show a central limit theorem. For that we must strengthen assumptions on $H$ that were demonstrated in the beginning of Section 4.2:

1. $H$ is of class $C^\infty$ and there exists a unique point $\theta_* \in [0, \pi[$ such that $H(\theta_*) = \underline{H} = \inf H$.
2. $H$ is decreasing on $]\theta_* - \pi/2, \theta_*[$ and $H$ is increasing on $]\theta_*, \theta_* + \pi/2[$. Moreover, $H''(\theta_*) > 0$.
3. $H$ is analytic in the neighborhood of the point $\theta_*$:

$$H(\theta) \overset{\theta \to \theta_*}{=} \underline{H} + \sum_{i=2}^{+\infty} H_i(\theta - \theta_*)^i.$$

4. The angle $\omega$ is such that $\omega \neq \theta_* + \pi/2 \pmod \pi$.

The assumptions imply that the function $\Lambda$ is also analytic in the neighborhood of the point $\theta_*$:

$$\Lambda(\theta) \overset{\theta \to \theta_*}{=} \sum_{i=0}^{+\infty} \Lambda_i(\theta - \theta_*)^i.$$

In this case, we use the extended Laplace method (see Wong [17]).

**Proposition 17.** *For all $q \in \mathbb{N}$, $q \geq 1$, we have, almost surely,*

$$n^{2\underline{H}-1}\sqrt{\log n}\, V_n(Z) \overset{n \to +\infty}{=} 16 \sum_{i=0}^{q} \frac{\Gamma((i+1)/2)\sigma_i}{(\log n)^{i/2}} t + o\left(\frac{1}{(\log n)^{q/2}}\right),$$

*where the coefficients $\sigma_i$ can be expressed in terms of $H_i$ and $\Lambda_i$.*

**Remark.** We can give that explicit forms for the coefficients $\sigma_i$ (see Wong [17]). For instance, the first two coefficients are given by

$$\sigma_0 = \frac{4 - 2^{2\underline{H}}}{16\sqrt{\pi}} G_{\theta_*}, \tag{90}$$

$$\sigma_1 = \frac{2}{H''(\theta_*)}\left(\frac{\Lambda'(\theta_*)(4 - 2^{2\underline{H}})}{2} - \frac{H^{(3)}(\theta_*)\Lambda(\theta_*)}{3H''(\theta_*)}\right). \tag{91}$$

**Proposition 18.** *We have, when $n \to +\infty$,*

$$\sqrt{n \log n}\left(n^{2\underline{H}-1} V_n(Z) - \int_0^\pi \Lambda(\theta)(4 - 2^{2H(\theta)})\left(\frac{1}{n}\right)^{2(H(\theta)-\underline{H})} \mathbf{1}_{\{H(\theta) \leq \underline{H} + 1/4\}}\, d\theta\right)$$

$$\overset{(\mathcal{L})}{\longrightarrow} \mathcal{N}(0, G^2_{\theta_*} \sigma^2_{\mathrm{FBM}, \underline{H}}). \tag{92}$$

Theorem 3 and the Slutsky lemma yield that when $n \to +\infty$,

$$\sqrt{n \log n}\left(\begin{array}{c} n^{2\underline{H}-1} V_n(Z) - \int_0^\pi \Lambda(\theta)(4 - 2^{2H(\theta)})(\frac{1}{n})^{2(H(\theta)-\underline{H})} \mathbf{1}_{\{H(\theta) \leq \underline{H} + 1/4\}}\, d\theta \\ (2n)^{2\underline{H}-1} V_{2n}(Z) - \int_0^\pi \Lambda(\theta)(4 - 2^{2H(\theta)})(\frac{1}{2n})^{2(H(\theta)-\underline{H})} \mathbf{1}_{\{H(\theta) \leq \underline{H} + 1/4\}}\, d\theta \end{array}\right)$$

$$\xrightarrow{(\mathcal{L})} \mathcal{N}(0, G_{\theta_*}^2 \Sigma_{\mathrm{FBM}, \underline{H}}). \tag{93}$$

As in the previous case, we obtain an estimator of $\underline{H}$ that is strongly consistent and asymptotically normal.

**Corollary 19.** *The statistic*

$$\widehat{\underline{H}}_n = \frac{1}{2} - \log \frac{V_{2n}(Z)}{V_n(Z)} (2\log 2)^{-1} \tag{94}$$

*is a strongly consistent estimator of $\underline{H}$ and*

$$\sqrt{n}\left(\widehat{\underline{H}}_n - \underline{H} + \log\left(\frac{G_{\theta_*} + \phi(1/(2N_n))}{G_{\theta_*} + \phi(1/N_n)}\right)(2\log 2)^{-1}\right)$$

$$\xrightarrow{(\mathcal{L})} \mathcal{N}\left(0, G_{\theta_*} \frac{3\sigma_{\mathrm{FBM}, \underline{H}}^2 - 2^{2-2\underline{H}}\sigma_{*, \mathrm{FBM}, \underline{H}}^2}{4(4 - 2^{2\underline{H}})\log 2}\right) \tag{95}$$

*with*

$$\phi(h) = 8\sqrt{-\log h} \int_0^\pi \Lambda(\theta)(4 - 2^{2H(\theta)})|h|^{2(H(\theta) - \underline{H})} \mathbf{1}_{\{H(\theta) \le \underline{H} + 1/4\}} \, \mathrm{d}\theta$$

$$- (4 - 2^{2\underline{H}}) G_{\theta_*}.$$

Let us note that because $Z$ is Gaussian and has stationary increments, we can apply the results of Istas and Lang [15] to estimate $\underline{H}$. To do so we must assume that $\mathrm{Leb}(\{H(\theta) = \underline{H}\}) > 0$ and $\underline{H} \ge 3/4$. Moreover, we need observations of $Z$ along an infinite interval, which is not the case in our assumptions ($t \in [0, 1]$). In this sense, we have improved the result of Istas and Lang [15] in the case of the AFBM.

Let us note that the estimation of the function $H$ was performed by Ayache *et al.* [2] and Biermé [6].

## 4.3. Proof of the results in the l.a.s.s. case

**Proof of Proposition 10.** We must check that the assumptions of Theorem 1 are satisfied. For assumption 2(a) (in Theorem 1), note that the functions $H$ and $\Lambda$ are bounded functions that can be deduced by the Lebesgue theorem of continuity under the symbol integral.

For assumption 2(b) (in Theorem 1), we must compute the derivative $\frac{\partial^4 R}{\partial s^2 \partial t^2}$. For the same reasons as above, the Lebesgue theorem of differentiability under the symbol integral shows the existence of this derivative on $]0, 1]^2 \setminus \{s = t\}$ and yields the formula

$$|s - t|^{2+\gamma} \frac{\partial^4 R}{\partial s^2 \partial t^2}(s, t)$$

$$= -4 \int_0^\pi \Lambda(\theta) 2H(\theta)(2H(\theta)-1)(2H(\theta)-2)(2H(\theta)-3)|s-t|^{2(H(\theta)-\underline{H})} \, \mathrm{d}\theta, \quad (96)$$

where $\gamma = 2(1 - \underline{H})$. Because the right-hand side is bounded, this proves that the assumption 2(b) (in Theorem 1) is satisfied with $L(h) = 1$.

Computing $(\delta_1^h \circ \delta_2^h R)(s,t)$, we get

$$(\delta_1^h \circ \delta_2^h R)(s,t)$$
$$= 4 \int_0^\pi \Lambda(\theta)(-|s-t-2h|^{2H(\theta)} + 4|s-t-h|^{2H(\theta)} - 6|s-t|^{2H(\theta)}$$
$$+ 4|s-t+h|^{2H(\theta)} - |s-t+2h|^{2H(\theta)}) \, \mathrm{d}\theta.$$

Thus

$$(\delta_1^h \circ \delta_2^h R)(t,t) = 8 \int_0^\pi \Lambda(\theta)|h|^{2H(\theta)}(4 - 2^{2H(\theta)}) \, \mathrm{d}\theta,$$

so

$$\frac{(\delta_1^h \circ \delta_2^h R)(t,t)}{|h|^{2-\gamma}} = 8 \int_0^\pi \Lambda(\theta)|h|^{2(H(\theta)-\underline{H})}(4 - 2^{2H(\theta)}) \, \mathrm{d}\theta \quad (97)$$

with $\gamma = 2 - 2\underline{H}$.

Setting $g_0(t) = (4 - 2^{2\underline{H}})J_{\underline{H}}$, we have

$$\frac{(\delta_1^h \circ \delta_2^h R)(t,t)}{|h|^{2-\gamma}} - g_0(t) = 8 \int_0^\pi \Lambda(\theta)\mathbf{1}_{\{H(\theta) > \underline{H}\}}(4 - 2^{2H(\theta)})|h|^{2(H(\theta)-\underline{H})} \, \mathrm{d}\theta. \quad (98)$$

Therefore, the Lebesgue theorem and the fact that the right-hand side does not depend on $t$ yield

$$\lim_{h \to 0^+} \sup_{h \le t \le 1-h} \left| \frac{(\delta_1^h \circ \delta_2^h R)(t,t)}{|h|^{2-\gamma}} - g_0(t) \right| = 0.$$

Hence the assumption 2(c) (in Theorem 1) is fulfilled.

So if we apply Theorem 1 to $Z$, we obtain (79). $\qquad \square$

**Proof of Proposition 11.** We apply Theorem 2 to $Z$. We must show that the assumptions 2 and 3 (in Theorem 2) are satisfied. For assumption 2 (in Theorem 2), we must compute the derivative $\frac{\partial^4 R}{\partial s^2 \partial t^2}$. As in the proof of Proposition 10, the Lebesgue theorem of differentiability under the symbol integral shows the existence of this derivative on $]0,1]^2 \setminus \{s = t\}$ and yields the formula (96),

$$|s-t|^{2+\gamma} \frac{\partial^4 R}{\partial s^2 \partial t^2}(s,t)$$
$$= -4 \int_0^\pi \Lambda(\theta) 2H(\theta)(2H(\theta)-1)(2H(\theta)-2)(2H(\theta)-3)|s-t|^{2(H(\theta)-\underline{H})} \, \mathrm{d}\theta,$$

where $\gamma = 2(1 - \underline{H})$. Hence the assumption 2 (in Theorem 2) is satisfied with $L(h) = 1$ and $C(s, t)$ equal to the right-hand side of (96) (the continuity of $C$ on $\{s = t\}$ is obtained by applying Lebesgue theorem).

Next we deal with the assumption 3 (in Theorem 2). We have proved (97):

$$\frac{(\delta_1^h \circ \delta_2^h R)(t, t)}{|h|^{2-\gamma}} = 8 \int_0^\pi \Lambda(\theta) |h|^{2(H(\theta) - \underline{H})} (4 - 2^{2H(\theta)}) \, \mathrm{d}\theta.$$

We set $q = 0$, $\alpha_0 = 1/2$ and

$$g_0(t) = 8(4 - 2^{2\underline{H}}) \int_0^\pi \Lambda(\theta) \mathbf{1}_{\{H(\theta) = \underline{H}\}} \, \mathrm{d}\theta = (4 - 2^{2\underline{H}}) J_{\underline{H}}.$$

We have

$$\frac{(\delta_1^h \circ \delta_2^h R)(t, t)}{|h|^{2-\gamma}} - g_0(t) = 8 \int_0^\pi \Lambda(\theta) |h|^{2(H(\theta) - \underline{H})} (4 - 2^{2H(\theta)}) \mathbf{1}_{\{H(\theta) > \underline{H} + 1/4\}} \, \mathrm{d}\theta.$$

Therefore, the Lebesgue theorem yields

$$\lim_{h \to 0^+} \frac{1}{\sqrt{h}} \sup_{h \le t \le 1 - h} \left| \frac{(\delta_1^h \circ \delta_2^h R)(t, t)}{|h|^{2-\gamma}} - g_0(t) \right| = 0$$

and assumption 3(d) (in Theorem 2) is fulfilled.

Moreover, we have

$$(\delta_1^h \circ \delta_2^h R)(t + h, t) = 4 \int_0^\pi \Lambda(\theta) |h|^{2H(\theta)} (4.2^{2H(\theta)} - 3^{2H(\theta)} - 7) \, \mathrm{d}\theta,$$

so

$$\frac{(\delta_1^h \circ \delta_2^h R)(t + h, t)}{|h|^{2-\gamma}} = 4 \int_0^\pi \Lambda(\theta) |h|^{2(H(\theta) - \underline{H})} (4.2^{2H(\theta)} - 3^{2H(\theta)} - 7) \, \mathrm{d}\theta.$$

Consequently, the theorem of dominated convergence yields

$$\lim_{h \to 0^+} \sup_{h \le t \le 1 - h} \left| \frac{(\delta_1^h \circ \delta_2^h R)(t + h, t)}{|h|^{2-\gamma}} - \widetilde{g}(t) \right| = 0,$$

where

$$\widetilde{g}(t) = \tfrac{1}{2}(2^{2\underline{H} + 2} - 3^{2\underline{H}} - 7) J_{\underline{H}}.$$

Thus, the assumption 3 (in Theorem 2) is fulfilled.                                                    $\square$

**Proof of Corollary 12.** The almost sure convergence is a straightforward consequence of (79). To get (82), we just apply the $\delta$ method with the $C^1$ function:

$$f(x, y) = -\frac{\log(y/x)}{2 \log 2}. \qquad\qquad\qquad \square$$

**Proof of Proposition 13.** We apply Theorem 1 to $Z$ to obtain a refinement of Proposition 10. We must show that Assumptions 2(b), 2(c) and 3 (in Theorem 1) are satisfied. For assumption 2(b) (in Theorem 1), we use the same arguments as in Proposition 11 and obtain the same function $C(s,t)$ with $\gamma = 2 - 2\underline{H}$ and $L(h) = 1$.

Next we deal with assumption 2(c) (in Theorem 1). Assumption 2(c)(iii) (in Theorem 1) is a straightforward consequence of the Lebesgue theorem.

Moreover, we have proved (97):

$$\frac{(\delta_1^h \circ \delta_2^h R)(t,t)}{|h|^{2-\gamma}} = 8 \int_0^\pi \Lambda(\theta) |h|^{2(H(\theta) - \underline{H})} (4 - 2^{2H(\theta)}) \, d\theta.$$

We set $q = 1$, $\alpha_1 = 1/2$, $\nu_1 = 1$ and

$$g_0(t) = 8(4 - 2^{2\underline{H}}) \int_0^\pi \Lambda(\theta) \mathbf{1}_{\{H(\theta) = \underline{H}\}} \, d\theta = J_{\underline{H}}(4 - 2^{2\underline{H}}),$$

$$g_1(t) = 1,$$

$$\phi(h) = 8 \int_0^\pi \Lambda(\theta)(4 - 2^{2H(\theta)}) |h|^{2(H(\theta) - \underline{H})} \mathbf{1}_{\{\underline{H} < H(\theta) \le \underline{H} + 1/4\}} \, d\theta.$$

We have

$$\frac{(\delta_1^h \circ \delta_2^h R)(t,t)}{|h|^{2-\gamma}} - g_0(t) - g_1(t)\phi(h)$$

$$= 8 \int_0^\pi \Lambda(\theta)(4 - 2^{2H(\theta)}) |h|^{2(H(\theta) - \underline{H})} \mathbf{1}_{\{\underline{H} + 1/4 < H(\theta)\}} \, d\theta$$

$$\overset{h \to 0^+}{=} o(\sqrt{h}),$$

thanks to Lebesgue theorem. Because the right-hand side does not depend on $t$, assumption 1(c)(v) (in Theorem 1) is fulfilled.

Moreover, we have

$$\phi(h) = 8 \int_0^\pi \Lambda(\theta)(4 - 2^{2H(\theta)}) |h|^{2(H(\theta) - \underline{H})} \mathbf{1}_{\{\underline{H} < H(\theta) \le \underline{H} + 1/4\}} \, d\theta.$$

Therefore, for $h$ enough small,

$$\phi(h) \ge K h^{2(\underline{H} + 1/4 - \underline{H})} = K\sqrt{h},$$

which yields that the assumption 2(c)(iii) (in Theorem 1) is satisfied too. The last inequality yields that the assumption 3 (in Theorem 1) is fulfilled.

Therefore, (83) is a consequence of Theorem 1 applied to $Z$. To prove (84), we apply Theorem 2. We need to check that assumptions 3(d) and 3(e) (in Theorem 2) are fulfilled. Assumption 3(d) (in Theorem 2) is a straightforward consequence of previous computations. For assumption 3(e) (in Theorem 2), we use the same arguments and the same function $\widetilde{g}(t)$ as in Proposition 11. □

**Proof of Corollary 14.** The almost sure convergence is a straightforward consequence of (79). To prove (86), we apply the $\delta$ method between the points $\binom{n^{2HK-1}V_n(Z)}{(2n)^{2HK-1}V_{2n}(Z)}$ and $\binom{(4-2^{2\underline{H}})J_{\underline{H}}+\phi(\frac{1}{n})}{(4-2^{2\underline{H}})J_{\underline{H}}+\phi(\frac{1}{2n})}$ to the $C^1$ function:

$$f(x,y) = -\frac{\log(y/x)}{2\log 2}.$$

□

## 4.4. Proofs of the results in the non-l.a.s.s. case

**Proof of Proposition 15.** Thanks to Proposition 9, we can compute the variogram of $A^{(H)}$ for all $t \in \mathbb{R}^2$:

$$v(t) = \frac{1}{8}\int_0^{2\pi} C(1, H(\theta))^2 |t_1\cos(\theta) + t_2\sin(\theta)|^{2H(\theta)}\,\mathrm{d}\theta.$$

We use the polar coordinates and we parametrize $t \in \mathbb{R}^2 \setminus \{(0,0)\}$ by $(\rho(t), \alpha(t))$:

$$t_1 = \rho(t)\cos\alpha(t),$$
$$t_2 = \rho(t)\sin\alpha(t).$$

We set $v_\varepsilon(t) = \varepsilon^{-2\underline{H}}v(\varepsilon t), \varepsilon > 0$, and use the polar parametrization and the $\pi$-periodicity of the function $H$ to obtain

$$v_\varepsilon(t) = \frac{1}{4}\int_{\theta_*-\pi/2}^{\theta_*+\pi/2} C(1, H(\theta))^2 \rho(t)^{2H(\theta)}|\cos(\alpha(t)-\theta)|^{2H(\theta)}\varepsilon^{2(H(\theta)-\underline{H})}\,\mathrm{d}\theta.$$

We assume that $\rho(t) \neq 0$ and $\alpha(t) \neq \theta_* + \pi/2 \pmod\pi$. The Laplace method (see Dieudonné [11], Theorem IV.2.5, page 125) yields

$$v_\varepsilon(t) \overset{\varepsilon\to 0^+}{\sim} \frac{C(1,\underline{H})^2}{8}\rho(t)^{2\underline{H}}|\cos(\alpha(t)-\theta_*)|^{2\underline{H}}\sqrt{\frac{\pi}{(-\log\varepsilon)H''(\theta_*)}}.$$

Therefore,

$$\lim_{\varepsilon\to 0^+}\sqrt{\log\left(\frac{1}{\varepsilon}\right)}\frac{1}{\varepsilon^{2\underline{H}}}v(\varepsilon t) = \frac{C(1,\underline{H})^2}{8}\rho(t)^{2\underline{H}}|\cos(\alpha(t)-\theta_*)|^{2\underline{H}}\sqrt{\frac{\pi}{H''(\theta_*)}}.$$

With a refinement of the Laplace method, we can check that this is still true if $\rho(t) = 0$ or $\alpha(t) = \theta_* + \pi/2 \pmod\pi$, which are the cases where the limit is vanishing.

Moreover, because $A^{(H)}$ is a Gaussian processes with stationary increments and vanishing limits at the origin, we have

$$\forall s, t \in \mathbb{R}^d \qquad \text{Cov}(A_s^{(H)}, A_t^{(H)}) = 2(v(s) + v(t) - v(s-t)).$$

We are able to conclude that the finite-dimensional laws of the process (87), with $t_0 = 0$, converge toward those of the process defined in (88). Because $A^{(H)}$ has stationary increments, we have the same result for all $t_0 \in \mathbb{R}^d$. □

**Proof of Proposition 16.** As in Proposition 10, we must check that the assumptions of Theorem 1 are satisfied. For assumption 2(b), in Proposition 10 we showed formula (96),

$$|s-t|^{2+\gamma} \frac{\partial^4 R}{\partial s^2 \partial t^2}(s, t)$$
$$= -4 \int_0^\pi \Lambda(\theta) 2H(\theta)(2H(\theta)-1)(2H(\theta)-2)(2H(\theta)-3)|s-t|^{2(H(\theta)-\underline{H})} \, \mathrm{d}\theta,$$

where $\gamma = 2(1 - \underline{H})$. We set $L(h) = \sqrt{\frac{-1}{\log h}}$, where $L$ is slowly varying, and write

$$\frac{|s-t|^{2+\gamma}}{L(|s-t|)} \frac{\partial^4 R}{\partial s^2 \partial t^2}(s, t)$$
$$= -4 \int_0^\pi \Lambda(\theta) 2H(\theta)(2H(\theta)-1)(2H(\theta)-2)(2H(\theta)-3)$$
$$\times \sqrt{\log(|s-t|)} |s-t|^{2(H(\theta)-\underline{H})} \, \mathrm{d}\theta. \tag{99}$$

If $\underline{H} \neq \frac{1}{2}$, the Laplace method yields that the right-hand side goes to

$$\sqrt{\frac{\pi}{H''(\theta_*)}} \Lambda(\theta_*) 2\underline{H}(2\underline{H}-1)(2\underline{H}-2)(2\underline{H}-3) \tag{100}$$

when $|s-t| \to 0$. A refinement of the Laplace method allows us to check that it is still true if $\underline{H} = \frac{1}{2}$. This implies that the right-hand side of (99) is a continuous function of $(s, t)$ on the set $[0, 1]^2$. Therefore, it is bounded and the assumption 2(b) (in Theorem 1) is satisfied.

For assumption 2(c) (in Theorem 1), we have seen in (97) that

$$\frac{(\delta_1^h \circ \delta_2^h R)(t, t)}{|h|^{2-\gamma}} = 8 \int_0^\pi \Lambda(\theta) |h|^{2(H(\theta)-\underline{H})} (4 - 2^{2H(\theta)}) \, \mathrm{d}\theta$$

with $\gamma = 2(1 - \underline{H})$. We want to study the asymptotic behavior of the preceding integral, so we denote by $I(h)$ the right-hand side of (97). The Laplace method yields

$$I(h) \overset{h \to 0^+}{\sim} 8\Lambda(\theta_*)(4 - 2^{2\underline{H}}) \sqrt{\frac{2\pi}{(-\log h) 2H''(\theta_*)}}.$$

Therefore, we obtain

$$\lim_{h \to 0^+} \sqrt{(-\log h)} I(h) = 8\Lambda(\theta_*)(4 - 2^{2\underline{H}})\sqrt{\frac{\pi}{H''(\theta_*)}}.$$

As in the Proof of Proposition 10, because $\sqrt{(-\log h)} I(h)$ does not depend on the variable $t$, we have

$$\lim_{h \to 0^+} \sup_{h \le t \le 1-h} \left| \frac{(\delta_1^h \circ \delta_2^h R)(t,t)}{L(h)|h|^{2-\gamma}} - g_0(t) \right| = 0,$$

where $g_0(t) = (4 - 2^{2\underline{H}})G_{\theta_*}$. Therefore, the assumption 2(c) (in Theorem 1) is fulfilled with $q = 0$. □

**Proof of Proposition 17.** We apply Theorem 1. The difference from the proof of Proposition 16 comes from assumption 2(c) (in Theorem 1). We must compute an asymptotic expansion of the expression, which is a consequence of (97),

$$\frac{(\delta_1^h \circ \delta_2^h R)(t,t)}{L(h)|h|^{2-\gamma}} = 4\sqrt{-\log h} \int_0^{2\pi} \Lambda(\theta)|h|^{2(H(\theta)-m)}(4 - 2^{2H(\theta)}) \, \mathrm{d}\theta$$

with $\gamma = 2(1 - \underline{H})$. The $\pi$-periodicity of the functions $H$ and $\Lambda$ yields

$$\frac{(\delta_1^h \circ \delta_2^h R)(t,t)}{L(h)|h|^{2-\gamma}} = 8\sqrt{-\log h} \int_{\theta_*}^{\theta_* + \pi} \Lambda(\theta)(4 - 2^{2H(\theta)})\mathrm{e}^{2(H(\theta)-m)\log h} \, \mathrm{d}\theta.$$

Here we cut this integral into two parts (we integrate on $[\theta_*, \theta_* + \pi/2]$ and on $[\theta_* + \pi/2, \theta_* + \pi]$) and use Theorem II.1.1 from Wong [17] on the extended Laplace method. We obtain

$$\frac{(\delta_1^h \circ \delta_2^h R)(t,t)}{L(h)|h|^{2-\gamma}} \overset{h \to 0^+}{=} 16 \sum_{i=0}^{q} \frac{\Gamma((i+1)/2)\sigma_i}{(-\log h)^{i/2}} + o\left(\frac{1}{(-\log h)^{q/2}}\right), \tag{101}$$

where the coefficients $\sigma_i$ can be expressed in terms of $H_i$ and $\Lambda_i$.

Because these quantities do not depend on the variable $t$, the assumption 2(c)(iii) of Theorem 1 is fulfilled with

$$g_i(t) = 16\Gamma\left(\frac{i+1}{2}\right)\sigma_i, \qquad \varepsilon_i = 1/2, \qquad \phi(h) = \frac{1}{\sqrt{-\log h}}, \qquad \nu_i = i. \qquad □$$

**Proof of Proposition 18.** We apply Theorem 2 to $Z$. We must show that assumptions 2 and 3 (in Theorem 2) are satisfied. As in Proposition 16, we can check that assumption 2 (in Theorem 2) is fulfilled with $C(s,t)$ equal to the right-hand side of (99).

For assumption 3 (in Theorem 2), we have proved (97):

$$\frac{(\delta_1^h \circ \delta_2^h R)(t,t)}{|h|^{2-\gamma}} = 8 \int_0^{\pi} \Lambda(\theta)|h|^{2(H(\theta)-\underline{H})}(4 - 2^{2H(\theta)}) \, \mathrm{d}\theta.$$

We set $q = 0$, $\alpha_1 = 1/2$, $\nu_1 = 1$ and

$$g_0(t) = (4 - 2^{2\underline{H}})G_{\theta_*},$$

$$g_1(t) = 1,$$

$$\phi(h) = 8\sqrt{-\log h} \int_0^\pi \Lambda(\theta)(4 - 2^{2H(\theta)})|h|^{2(H(\theta)-\underline{H})} \mathbf{1}_{\{H(\theta) \leq \underline{H}+1/4\}} \, d\theta$$

$$- (4 - 2^{2\underline{H}})G_{\theta_*}.$$

We have

$$\frac{(\delta_1^h \circ \delta_2^h R)(t,t)}{L(h)|h|^{2-\gamma}} - g_0(t) - g_1(t)\phi(h)$$

$$= 8\sqrt{-\log h} \int_0^\pi \Lambda(\theta)(4 - 2^{2H(\theta)})|h|^{2(H(\theta)-\underline{H})} \mathbf{1}_{\{\underline{H}+1/4 < H(\theta)\}} \, d\theta$$

$$\stackrel{h \to 0^+}{=} o(\sqrt{h}),$$

thanks to the theorem of dominated convergence. Because the right-hand side does not depend on $t$, assumption 3(d) (in Theorem 2) is fulfilled.

Moreover, we have

$$(\delta_1^h \circ \delta_2^h R)(t+h, t) = 4 \int_0^\pi \Lambda(\theta)|h|^{2H(\theta)}(4.2^{2H(\theta)} - 3^{2H(\theta)} - 7) \, d\theta$$

and so

$$\frac{(\delta_1^h \circ \delta_2^h R)(t+h, t)}{L(h)|h|^{2-\gamma}} = 4\sqrt{-\log h} \int_0^\pi \Lambda(\theta)|h|^{2(H(\theta)-\underline{H})}(4.2^{2H(\theta)} - 3^{2H(\theta)} - 7) \, d\theta.$$

Consequently, the Laplace method and the theorem of dominated convergence yield

$$\lim_{h \to 0^+} \sup_{h \leq t \leq 1-h} \left| \frac{(\delta_1^h \circ \delta_2^h R)(t+h, t)}{L(h)|h|^{2-\gamma}} - 4\Lambda(\theta_*)(4.2^{2\underline{H}} - 3^{2\underline{H}} - 7)\sqrt{\frac{\pi}{H''(\theta_*)}} \right| = 0.$$

Thus, the assumption 3(e) (in Theorem 2) is fulfilled with

$$\widetilde{g}(t) = \frac{4.2^{2\underline{H}} - 3^{2\underline{H}} - 7}{2} G_{\theta_*}. \qquad \square$$